\documentclass[11pt]{amsbook}

\usepackage{amsmath,amssymb,epsf, epic, eepic,fancybox}
\usepackage{psfig}
\usepackage{amsfonts}
\usepackage{times}
\usepackage{slashbox}

\input xy
\xyoption{all}

\newtheorem{thm}{Theorem}[section]

\newtheorem{prop}[thm]{Proposition}
\newtheorem{exe}[thm]{Example}

\theoremstyle{definition}

\newtheorem{exercise}[thm]{Exercise}

\newcommand{\bZ}{\mathbb{Z}}
\newcommand{\bQ}{\mathbb{Q}}
\newcommand{\bR}{\mathbb{R}}

\newcommand{\cC}{\mathcal{C}}

\newcommand{\cM}{\mathcal{M}}

\DeclareMathOperator{\Tor}{Tor}

\newcommand{\ra}{\rightarrow}
\newcommand{\ol}[1]{\overline{#1}}

\newcommand{\spaces}{\;\;\;\;\;\;\;}
\newcommand{\ot}{\otimes}
\newcommand{\lk}{\text{lk}}

\newcommand{\inlinepic}[1]{\raisebox{-1mm}{\psfig{figure={#1}}}}
\newcommand{\textpic}[1]{\raisebox{-3mm}{\psfig{figure={#1}}}}
\newcommand{\textpicc}[1]{\raisebox{-4mm}{\psfig{figure={#1}}}}
\newcommand{\inlinebigpic}[1]{\raisebox{-7mm}{\psfig{figure={#1}}}}

\newcommand{\cc}[3]{C^{{#1},{#2}}(#3)}
\newcommand{\ccc}[4]{C^{{#1},{#2}}(#3;#4)}
\newcommand{\kh}[3]{K\!H^{{#1},{#2}}(#3)}

\newcommand{\khc}[4]{K\!H^{{#1},{#2}}(#3;#4)}
\newcommand{\lee}[2]{Lee^{{#1}}(#2)}
\newcommand{\dimn}[1]{\text{dim}({#1})}
\newcommand{\qdim}[1]{\text{qdim}({#1})}
\newcommand{\cob}{\cC ob_{1+1}}
\newcommand{\Rmod}{\cM od_R}
\newcommand{\cano}{\frak{s}_\theta}
\newcommand{\canoo}{\frak{s}_{\theta_0}}
\newcommand{\canooo}{\frak{s}_{\theta_1}}
\newcommand{\smin}{s_{\text{min}}(K)}
\newcommand{\smax}{s_{\text{max}}(K)}
\newcommand{\smaxu}{s_{\text{max}}(U)}

\title{Five Lectures on Khovanov Homology}
\author{Paul Turner}

\begin{document}


\maketitle

\vspace*{0.1cm}

\begin{abstract}
These five lectures provide an introduction to Khovanov homology covering the basic definitions, important properties, a number of
variants and some applications. At the end of each lecture the reader
is referred to the relevant literature for further reading.
\end{abstract}

\vspace*{2.5cm}

\begin{center}
{\bf About these lectures}
\end{center}
These lectures were designed for the summer school {\em Heegaard-Floer
homology and Khovanov homology} in Marseilles, 29th May - 2nd June, 2006.

The intended audience is graduate students with some minimal
background in low-dimensional and algebraic topology.  I hesitated to
produce lecture notes at all, since much of the literature in the
subject is very well written, but decided in the end that notes could
serve some purpose.

In order to keep the narrative flowing I found it convenient to delay
all attributions of credit until the end of each lecture. I have
attempted to do this as accurately as possible and if I have
failed to properly attribute a certain piece of work or omitted to
mention someone in a particular context, my apologies to the injured
party in advance.

At the present time the pace of development of the subject is very
rapid and the reader is encouraged to consult math/GT for the latest
developments.

Many thanks F. Costantino, M. Mackaay and P. Vaz for their comments on a draft
version and to D. Matignon for organising a very stimulating summer school.

\vfill

\hfill Paul Turner \hspace*{32ex}

\hfill School of Mathematical and Computer Sciences 

\hfill Heriot-Watt
  University \hspace*{22ex}

\hfill Edinburgh EH14 4AS \hspace*{23ex}

\hfill Scotland \hspace*{35ex}

\vspace*{10mm}
\noindent
\hfill \email{paul@ma.hw.ac.uk}\hspace*{26ex}

\tableofcontents


\chapter{Lecture One}
In this lecture we begin with a very brief introduction to the
subject, followed by some recollections about the Jones polynomial. We
then define the main object of interest: the Khovanov complex of an
oriented link diagram.

\section{What is it all about?}
Given an oriented link {\em diagram}, $D$, Khovanov constructs in a purely combinatorial way a bi-graded chain complex $\cc **D$ associated to $D$.
\[
\xymatrix{D  \ar@{~>}[rr]^{\text{Khovanov}} && \cc **D}
\]
Given a chain complex we can apply {\em homology} to it and for $\cc **D$ this results in the {\em Khovanov homology}, $\kh **D$, of the diagram $D$. 
\[
\xymatrix{\cc **D  \ar@{~>}[rr]^{\text{Homology}} && \kh **D
}
\]
The following properties are satisfied:
\begin{enumerate}
\item If $D$ is related to another diagram $D^\prime$ by a sequence of Reidemeister moves then there is an isomorphism
\[
\kh **D \cong \kh ** {D^\prime}.
\]
\item The {\em graded Euler characteristic} is the unnormalised Jones polynomial i.e.
\[
\sum_{i,j\in \bZ} (-1)^iq^j \dimn{\kh ijD} = \hat{J}(D).
\]
\end{enumerate}

You should think, by way of analogy, of the relationship of the
ordinary Euler characteristic to homology. For a space $M$ (there
are some restrictions on $M$, say, a finite CW complex) we can assign a
numerical invariant, the Euler characteristic $\chi(M)\in\bZ$. This
can be calculated by a simple algorithm based on some combinatorial
information about the space (eg the CW structure, a triangulation
etc). On the other hand, homology assigns a graded vector space to $M$ and is related to the Euler characteristic by the formula:
\[
\sum (-1)^i \dimn{H_i(M;\bQ)} = \chi(M).
\]
In this way homology {\em categorifies} the Euler characteristic: a number
gets replaced by a (graded) vector space whose (graded) dimension
gives back the number you started with.

In the case of links we can assign quantum invariants such as the Jones
polynomial. These can be calculated by a simple algorithm based on some
combinatorial information about the link (e.g. a diagram). Khovanov
homology is the analogue of homology and categorifies the Jones
polynomial.

Homology has many advantages over the Euler characteristic. For example
\begin{itemize}
\item homology is a stronger invariant than the Euler characteristic,
\item homology reveals richer information e.g. torsion,
\item homology is a {\em functor}.
\end{itemize}
As we will see soon, Khovanov homology has similar advantages over the Jones polynomial.

\section{Recollections about the Jones polynomial}
Let $L$ be an oriented link and $D$ a diagram for $L$ with $n$
crossings. Suppose that $n_-$ of these are negative crossings (like
this: \inlinepic{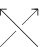}) and $n_+$ are positive (like this:
\inlinepic{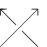}). The {\em Kauffman bracket} of the diagram $D$, written $\langle D
\rangle$ is the Laurent polynomial in a variable $q$ (i.e. $\langle D \rangle \in \bZ[q^{\pm 1}]$) defined recursively by:
\begin{eqnarray}
\langle \inlinepic{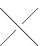} \rangle & = & \langle \inlinepic{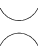} \rangle - q \langle \inlinepic{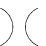} \rangle \label{eq:kb1}\\
\langle k \text{ circles in the plane } \rangle &= & (q+q^{-1})^k \label{eq:kb2}
\end{eqnarray}
(Beware: this is not the usual normalisation).

The Kauffman bracket is not a link invariant but by defining
\[
\hat{J}(D) = (-1)^{n_-}q^{n_+-2n_-}\langle D \rangle
\]
one gets a genuine link invariant i.e. something invariant under all the Reidemeister moves. The {\em Jones polynomial} is given by 
\[
J(D)= \frac{\hat{J}(D)}{q+q^{-1}}.
\]
(The usual formula for the Jones polynomial involves a variable $t$. Substitute $q = -t^{\frac{1}{2}}$ to make the descriptions match). The polynomial $\hat{J}(D)$ is known as the {\em unnormalised Jones polynomial}.

Equation (\ref{eq:kb1}) reduces the number of crossings at the expense
of twice as many terms on the right hand side. For a diagram $D$ with
$n$ crossings we can apply this equation $n$ times to end up with
$2^n$ pictures on the right hand side each of them consisting of a
collection of circles in the plane which we then evaluate using
Equation (\ref{eq:kb2}).

To do this in a systematic way let us agree that given a crossing
(looking like this: \inlinepic{crossing.eps}) we will call the two
pictures on the right of Equation (\ref{eq:kb1}) the {\em 0-smoothing}
(looking like this:
\inlinepic{smoothing0.eps}) and the {\em 1-smoothing} (looking like this:
\inlinepic{smoothing1.eps}).

\begin{figure}[h]
\centerline{
\psfig{figure= 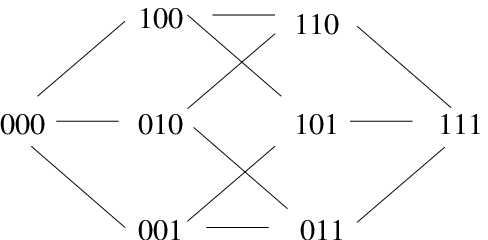}
}
\caption{}
\label{fig:cube1} 
\end{figure}

Thus, if we number the crossings of $D$ by $1,2, \ldots, n$ then each of the $2^n$ pictures can be indexed by a word of $n$ zeroes and ones i.e. an element of $\{0,1\}^n$. We will call a picture in which each crossing has been resolved (in one of the two ways above) a {\em smoothing}. Thus a diagram $D$ has $2^n$ smoothings indexed by $\{0,1\}^n$. The set $\{0,1\}^n$ is the vertex set of a hyper-cube as shown in Figure \ref{fig:cube1}, with an edge between words differing in exactly one place. 

It is convenient to think of the smoothings as living on the vertices
of this cube. For $\alpha\in\{0,1\}^n$ we will denote the associated
smoothing (the collection of circles in the plane) by $\Gamma_\alpha$.
Given $\alpha\in \{0,1\}^n$ we define
\[
r_\alpha = \text{ the number of 1's in }\alpha
\]
and 
\[
k_\alpha = \text{ the number of circles in }\Gamma_\alpha.
\]
We can now use Equations (\ref{eq:kb1}) and (\ref{eq:kb2}) to write down a state-sum expression for $\hat{J}(D)$.
\[
\hat{J}(D) = \sum_{\alpha\in \{0,1\}^n}(-1)^{r_\alpha +n_-}q^{r_\alpha + n_+ - 2n_-} (q+q^{-1})^{k_\alpha}.
\]

\begin{exercise}
Convince yourself that this above state-sum formula is correct.
\end{exercise}

\begin{exe}
Consider the Hopf link \inlinepic{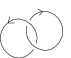} which has $n=n_-=2$ and $n_+=0$. There are four pictures which we assemble on the vertices of a square as displayed in Figure \ref{fig:hopflinkcube} and we compute the unnormalised Jones polynomial to be
\begin{eqnarray*}
\hat{J}(\inlinepic{hopflink.eps}) & = & q^{-4}(q+q^{-1})^2 -2 q^{-3}(q+q^{-1}) +q^{-2}(q+q^{-1})^2\\
& = & q^{-6} + q^{-4} + q^{-2} + 1.
\end{eqnarray*}

\begin{figure}[h]
\centerline{
\psfig{figure= 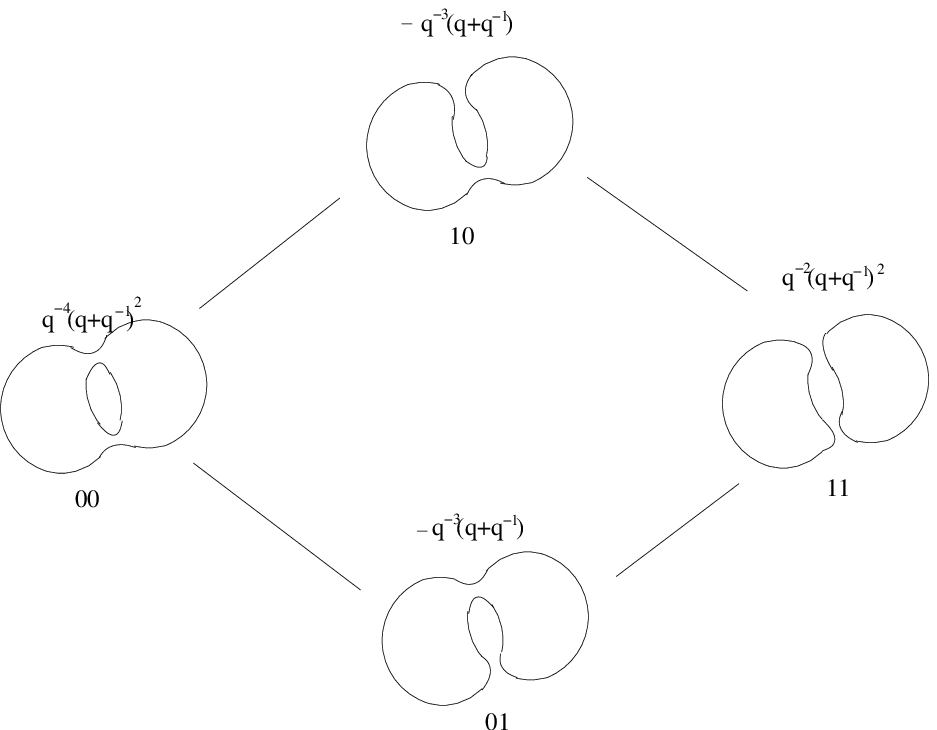}
}
\caption{}
\label{fig:hopflinkcube} 
\end{figure}

\end{exe}

\begin{exercise}
Write out the cube of smoothings for \inlinepic{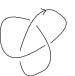} and use the state-sum formula to compute the unnormalised Jones polynomial.
\end{exercise}

\section{The definition of the Khovanov complex of a link diagram}
For the time being we will work over $\bQ$, so ``vector space'' means
`` vector space over $\bQ$''. Khovanov's idea is to assign a {\em cochain
complex} $(\cc **D,d)$ to a link diagram $D$. The homotopy type of
this complex will turn out to be an invariant and its graded Euler
characteristic the unnormalised Jones polynomial. 

Before getting to the definition let's recall a few things about finite dimensional graded vector spaces.

\begin{enumerate}
\item The {\em graded (or quantum) dimension}, qdim,  of a graded vector space $W=\bigoplus_m W^m$ is the polynomial in $q$ defined by
\[
\qdim W = \sum_m q^m \dimn {W^m}.
\]
\item The graded dimension satisfies
\[
\qdim {W\otimes W^\prime} = \qdim W \qdim {W^\prime},
\]
\[
\qdim {W\oplus W^\prime} = \qdim W + \qdim {W^\prime}.
\]
\item For a graded vector space $W$ and an integer $l$ we can define a new graded vector space $W\{l\}$ (a shifted version of $W$) by
\[
W\{l\}^m = W^{m-l}
\]
Notice that $\qdim {W\{l\}} = q^l\qdim W$.
\end{enumerate}

Now we turn to the definition of the {\em Khovanov complex}, $\cc
**D$, of an oriented link diagram $D$. An important role is played by
the following two-dimensional graded vector space. Let
$V=\bQ\{1,x\}$ (the $\bQ$-vector space with basis $1$ and $x$) and grade the two basis elements by
$\text{deg}(1) = 1$ and $\text{deg}(x)=-1$. 

\begin{exercise}
Show that $\qdim{V^{\otimes k}}=(q+q^{-1})^k$.
\end{exercise}

Recall that we have $2^n$
smoothings of our diagram. To each $\alpha\in\{0,1\}^n$ now
associate the graded vector space
\[
V_\alpha = V^{\otimes k_\alpha}\{r_\alpha +n_+ - 2n_-\}
\]
and define
\[
\cc i*D = \bigoplus_{\substack{\alpha\in\{0,1\}^n\\r_\alpha = i+n_-}}V_\alpha.
\]
The internal grading comes from the fact that each $V_\alpha$ is a
graded vector space. Note that the vector spaces $\cc i * D$ are
trivial outside the range $i=-n_-,\ldots, n_+$.

\begin{figure}[h]
\centerline{
\psfig{figure=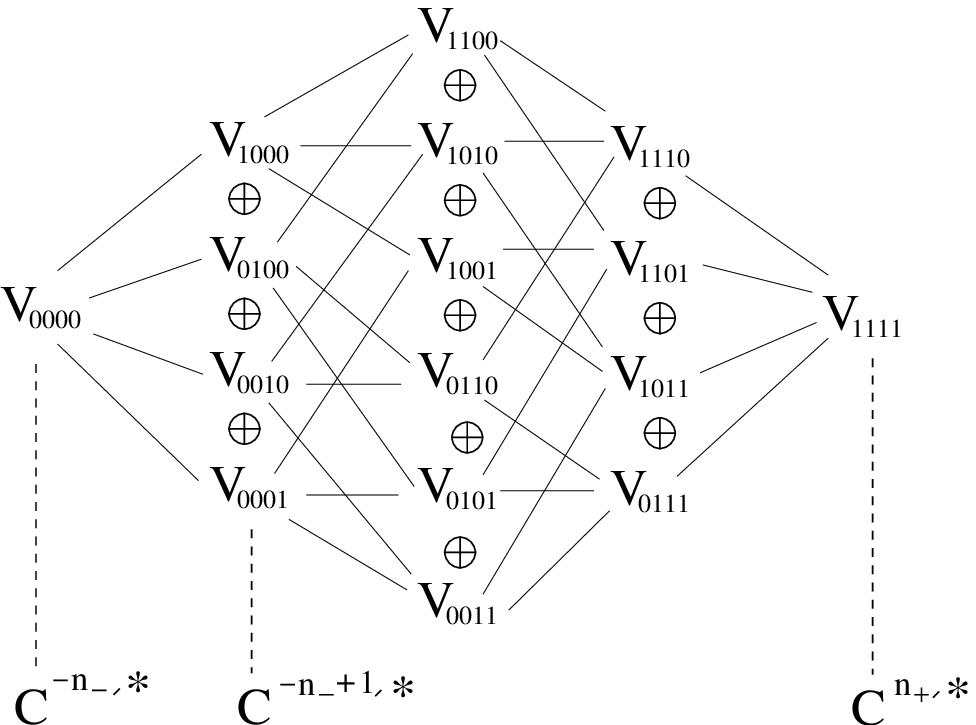}
}
\caption{}
\label{fig:cube2} 
\end{figure}

Recall that in the last section we arranged the $2^n$ smoothings of
the diagram on a cube with $2^n$ vertices indexed by $\{0,1\}^n$. The
above definition means that we now replace the smoothing
$\Gamma_\alpha$ with the vector space $V_\alpha$ so that the space
$\cc i * D$ is the direct sum of vector spaces in column $i+n_-$ of
the cube as indicated in Figure \ref{fig:cube2}.

\begin{exe}
For the Hopf link we have the cube in Figure \ref{fig:hopflinkcube2}.

\begin{figure}[h]
\centerline{
\psfig{figure= 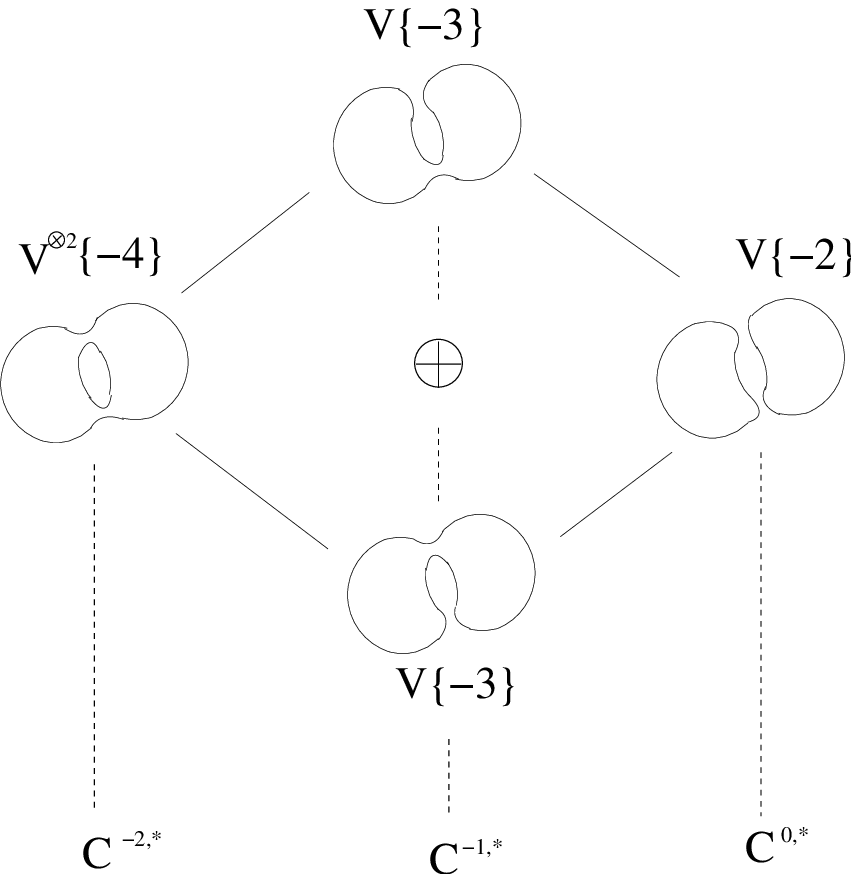}
}
\caption{}
\label{fig:hopflinkcube2} 
\end{figure}
\end{exe}

An element of $\cc ijD$is said to have {\em homological grading} $i$
and {\em $q$-grading} $j$. If $v\in V_\alpha\subset \cc**D$ with
homological grading $i$ and $q$-grading $j$ then it is useful to
remember that
\begin{eqnarray*}
i & = & r_\alpha - n_-\\
j & = & \text{deg}(v) +i +n_+ - n_-
\end{eqnarray*}
where deg$(v)$ is the degree of $v$ as an element of $V_\alpha$.

What we need now is a differential $d$ turning $(\cc **D, d)$ into a
complex.  Recall that we have a smoothing $\Gamma_\alpha$ (i.e. a
collection of circles) associated to each vertex $\alpha$ of the cube
$\{0,1\}^n$. Now to each edge of the cube we associate a {\em
cobordism} (i.e. an (orientable) surface whose boundary is the union of
the circles in the smoothings at either end).

\begin{figure}[h]
\centerline{
\psfig{figure= 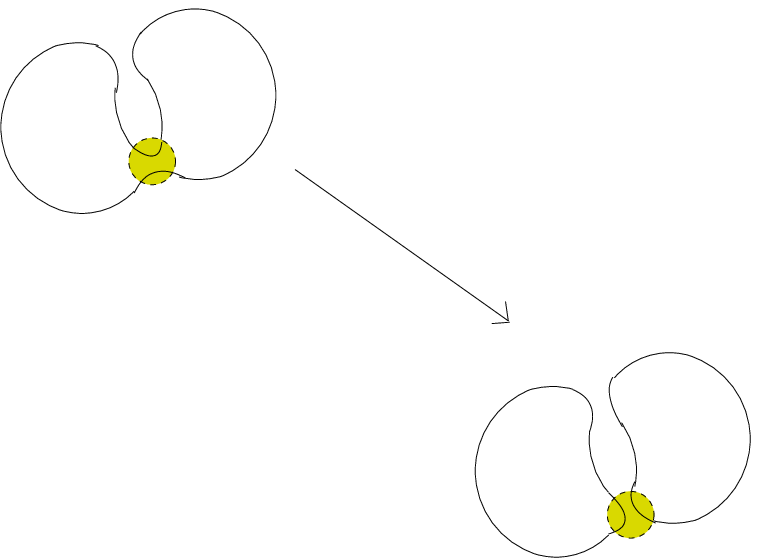}
}
\caption{}
\label{fig:chaningdiscHL} 
\end{figure}

Edges of the cube can be labelled by a string of zeroes and ones with
a star ($\star$) at the position that changes. For example the edge
joining $0100$ to $0110$ is denoted $01\star 0$. We can turn edges into {\em
arrows} by the rule: $\star = 0$ gives the tail and $\star=1$ gives
the head. For an arrow $\xymatrix{\alpha \ar[r]^\zeta &
\alpha^\prime}$ note that the smoothings $\alpha$ and $\alpha^\prime$
are identical except for a small disc (the {\em changing disc}) around
the crossing that changes from a 0- to a 1-smoothing (the one marked
by a $\star$ in $\zeta$). For example the changing disc for the arrow
$\zeta=1\star$ in the cube (square!) of the Hopf link above is shown
in Figure \ref{fig:chaningdiscHL}.

The cobordism $W_\zeta$ associated to $\xymatrix{\alpha \ar[r]^\zeta &
\alpha^\prime}$ is defined to be the following surface: outside the changing disc take the product of $\Gamma_\alpha$ with the unit interval and then plug the missing tube with  the saddle \textpic{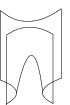}
 Thus each $W_\zeta$ consists
 of a bunch of cylinders and one pair-of-pants surface
 (\textpic{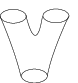} or
\textpic{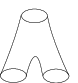}).

\vspace*{5mm}
{\bf Cobordism convention:} pictures of cobordisms go {\em down} the page.
\vspace*{5mm}

Above we replaced the smoothing $\Gamma_\alpha$ by the vector space
$V_\alpha$ and now we will replace the cobordism $W_\zeta$ associated to the edge $\xymatrix{\alpha
\ar[r]^\zeta & \alpha^\prime}$ by a linear
map $d_\zeta\colon V_\alpha \ra V_{\alpha^\prime}$. Since each circle
in a smoothing has a copy of the vector space $V$ attached to it, to
define $d_\zeta$ we only require two linear maps: one that fuses
$m\colon V\otimes V \ra V$ and one that splits $\Delta\colon V \ra
V\otimes V$. Then we can define $d_\zeta$ to be the identity on
circles not entering the changing disc and either $m$ or $\Delta$ on
the circles appearing in the changing disc (depending on whether the
pair-of-pants has two or one input boundary circles).

We define  $m\colon V\otimes V \ra V$ by
\[
1^2=1, \spaces 1x=x1=x, \spaces x^2=0,
\]
and $\Delta\colon V \ra V\otimes V$ by
\[
\Delta(1) = 1 \otimes x + x\otimes 1, \spaces \Delta(x) = x\otimes x.
\]

In fact by defining a unit $i(1) = 1$ and counit $\epsilon(1) = 0$ and
$\epsilon(x) = 1$ we have endowed $V$ with the structure of a
commutative {\em Frobenius algebra}. Isomorphism classes of commutative
Frobenius algebras are in bijective correspondence with isomorphisms
classes of 1+1-dimensional topological quantum field theories so what we
are really doing here is applying a TQFT (the one defined by $V$) to
the cube of circles and cobordisms - more on this in Lecture 2.

We are finally ready to define $d^i\colon \cc i*D \ra \cc {i+1}*D$. For $v\in V_\alpha\subset \cc i*D$ set
\[
d^i(v) = \sum_{\substack{\zeta \text{ such that}\\ \text{Tail} (\zeta)=\alpha}} \text{sign}(\zeta) d_\zeta(v)
\]
where $\text{sign}(\zeta)= (-1)^{\text{number of 1's to the left of }\star \text{ in } \zeta}$.

\begin{prop}
$d^{i+1}\circ d^i=0$.
\end{prop}
\begin{proof}(sketch)
The idea of the proof is that without the signs each face of the cube
commutes. To see this one can either look at a number of
cases and use the definition of the maps $m$ and $\Delta$ or (much
better) begin to think geometrically: each of the two routes around a
face gives the same cobordism (up to homeomorphism) and so applying
the TQFT defined by $V$ gives the same linear map.
Once all the non-signed faces commute then observe that the signs
occur in odd numbers on every face, thus turning commutativity into
anti-commutativity.
\end{proof}

\begin{exercise}
Write out the above proof properly - you may wish to wait until after
Lecture 2 where there is a more detailed discussion of Frobenius
algebras and TQFTs.
\end{exercise}

\begin{exercise}
Check that $d$ has bi-grading $(1,0)$.
\end{exercise}

The graded Euler characteristic of this complex i.e. 
\[
\sum (-1)^i \qdim {\cc i*D} \in \bQ[q^{\pm 1}]
\]
is nothing other than the unnormalised Jones polynomial.

\begin{exercise}
Convince yourself that the previous statement is true - this is simply
a matter of unwinding the definitions and then comparing with the
state-sum formula for the unnormalised Jones polynomial.
\end{exercise}

Later we will see that the homotopy type of $(\cc **D, d)$ is
invariant under transformations by Reidemeister moves. For now,
to end this lecture, let us perform a homology calculation.

\begin{exe}
Let us compute the {\em homology} of $\cc **
{\inlinepic{hopflink.eps}}$. The complex has only three non-trivial terms:
\[
\xymatrix{
0\ar[r] & \cc{-2}*D \ar[r]^d & \cc {-1}* D \ar[r]^d & \cc 0*D \ar[r]&  0.
}
\]
More explicitly we have:
\[
\xymatrix{
 & V\{-3\} \ar[dr]^{-\Delta} &\\
(V\otimes V)\{-4\} \ar[ur]^m \ar[dr]_m & \oplus & (V\otimes V)\{-2\}\\
&  V\{-3\} \ar[ur]^{\Delta} &
}
\]
and based on this one can compute as follows.

\vspace{5mm}

\begin{tabular}{|l|c|c|c|}
\hline
Homological degree & -2 & -1 & 0\\
\hline
Cycles & $\{1\otimes x - x\otimes 1,x\otimes x\}$ &
 $\{(1,1),(x,x)\}$ & $\{1\otimes 1,1\otimes x,x \otimes 1, x\otimes x\}$ \\
\hline
Boundaries & - & $\{(1,1),(x,x)\}$ & $\{1\otimes x + x\otimes 1, x\otimes x\}$\\
\hline
Homology &  $\{1\otimes x - x\otimes 1,x\otimes x\}$ & - & $\{1\otimes 1,1\otimes x\}$\\
\hline
$q$-degrees & -4,\hspace{1cm} -6 & & 0,\hspace{0.4cm} -2\\
\hline
\end{tabular}

\vspace{5mm}

The homology can be summarised as a table where the homological degree is horizontal and the $q$-degree vertical.

\begin{center}
\begin{tabular}{|c||c|c|c|}
\hline
\backslashbox{$j$}{$i$}& -2 & -1& 0\\
\hline\hline
0 & && $\bQ$\\\hline
-1 &&&\\\hline
-2 &&& $\bQ$\\\hline
-3&&&\\\hline
-4&$\bQ$&&\\\hline
-5&&&\\\hline
-6&$\bQ$&&\\\hline
\end{tabular}
\end{center}
\end{exe}

\begin{exercise}
Write out the cube for the trefoil explicitly (including the correct
signs on the edges). Calculate the homology of the Khovanov
complex. This involves a bit of work, but is a very good test to see
if you have understood all the definitions.
\end{exercise}

\section{Notes and further reading}
\noindent
The original paper by Khovanov in which he defines the complex and
establishes the basic properties is \cite{khovanov1}. In this paper he
starts out working over the ring $\bZ[c]$ and then sets $c=0$ to work
over $\bZ$. We work over $\bQ$ because some things are a little
simpler. We will say more about other coefficients in Lecture 3.
Bar-Natan's exposition of Khovanov's work \cite{barnatan2} is
extremely readable and has also been very influential.


\chapter{Lecture Two}
The complex $\cc **D$ for a link diagram $D$ defined in Lecture 1
depends very much on the diagram. However, it turns out that different
diagrams for the same link give complexes which are {\em homotopy
equivalent}. In this lecture we begin with an aside on Frobenius
algebras and topological quantum field theories and after this discuss
the homotopy invariance properties of the Khovanov complex
concentrating on the first Reidemeister move. Finally in this lecture
we define Khovanov homology and discuss some properties.

\section{Frobenius algebras and TQFTs}

Hidden in the background in Lecture 1, about to come to deserved
prominence, are 1+1-dimensional topological quantum field theories and
their algebraic counterparts, Frobenius algebras.

A {\em commutative Frobenius algebra over $R$} (a commutative ring
with unit) is a unital, commutative $R$-algebra $V$ which as an
$R$-module is projective of finite type (if $R=\bQ$ this just means a
finite dimensional vector space over $\bQ$), together with a module
homomorphism, the counit, $\epsilon\colon V \ra R$ such that the bilinear form
$\langle - , - \rangle\colon V\otimes V \ra R$ defined by $\langle v,
w \rangle = \epsilon (vw)$ is non-degenerate i.e. the adjoint
homomorphism $V\ra V^*$ is an isomorphism. It is useful to define a coproduct
$\Delta\colon V \ra V\ot V$ by $\Delta (v) = \sum_i v_i^\prime \ot
v_i^{\prime\prime}$ being the unique element such that for all $w\in
V$, $vw = \sum_i v_i^\prime \langle v_i^{\prime\prime},w\rangle$.

Frobenius algebras reflect the topology of surfaces. This statement is the rough equivalent of the more accurate:
\[
\{\text{Iso. classes of comm. Frobenius algebras}\} \longleftrightarrow 
\{\text{Iso. classes of 1+1-dimensional TQFTs}\}
\]
Recall that a 1+1-dimensional TQFT is a monoidal functor $\cob \ra
\Rmod$ where $\cob$ is the category whose objects are closed, oriented
1-manifolds and where a morphism $\Gamma \ra \Gamma^\prime$ is an
oriented surface $W$ with $\partial W = \ol{\Gamma} \sqcup
\Gamma^\prime$ (here the overline means take the opposite
orientation). In fact you have to be careful to get the
details of all this right - see the references at the end of the
lecture. The upshot is that a TQFT:
\begin{itemize}
\item assigns to each closed 1-manifold $\Gamma$, an $R$-module
 $V_\Gamma$ such that if $\Gamma = \Gamma_0 \sqcup \Gamma_1$ then $V_\Gamma
 = V_{\Gamma_0} \otimes V_{\Gamma_1}$ (this what the adjective
 ``monoidal'' refers to) and
\item assigns to each cobordism $W\colon \Gamma \ra \Gamma^\prime$, an $R$-homomorphism $V_\Gamma \ra V_{\Gamma^\prime}$.
\end{itemize}
These assignments are subject to some axioms which, among other things, guarantee
\begin{itemize}
\item homeomorphic cobordisms induce the same homomorphism,
\item gluing of cobordisms is well behaved, and
\item $V_\emptyset = R$.
\end{itemize}
 
\begin{figure}[h]
\centerline{
\psfig{figure= 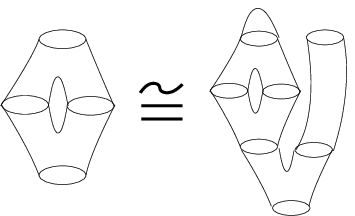}
}
\caption{}
\label{fig:s1} 
\end{figure}

Using the correspondence between Frobenius algebras and
1+1-dimensional TQFTs one can use the topology to prove algebraic
statements. For example, for any Frobenius algebra one can check that
$m(\Delta(v)) = m(m(\Delta(1)),v)$ for all $v\in V$. Checking this
algebraically is a bit of a pain, but geometrically it is a
triviality: the two surfaces in Figure \ref{fig:s1} are homeomorphic
and therefore correspond to the same homomorphism of Frobenius algebras.

In Lecture 1 we defined a particular two dimensional Frobenius algebra
$V$ which, by the above, defines a 1+1-dimensional TQFT. Given a
link diagram we considered the cube $\{0,1\}^n$ and associated to each
vertex $\alpha$ a collection of circles $\Gamma_\alpha$ (a smoothing) and to each edge $\zeta$ a
cobordism $W_\zeta$. In order to get a complex we then replaced each collection of circles by a vector
space and each cobordism by a linear map. This last step is nothing
other than applying the TQFT defined by $V$.

\begin{exercise}
A TQFT associates to the empty manifold the ground ring $R$ and thus a
closed cobordisms gives an element of $R$ (a closed cobordism is a
cobordism $\emptyset$ to $\emptyset$ and hence induces a map $R\ra R$
which you evaluate at $1\in R$). Compute the value of the torus for
the TQFT defined by the Frobenius algebra $V$ of Lecture 1.
\end{exercise}

\section{Reidemeister invariance}
 Recall that complexes $A^*$ and $B^*$ are homotopy
equivalent if there are chain maps $F\colon A^*\ra B^*$ and $G\colon
B^*\ra A^*$ such that $GF-Id_{A^*}$ and $FG-Id_{B^*}$ are
null-homotopic. 

\begin{prop}
If $D^\prime$ is a diagram obtained from $D$ by
the application of a Reidemeister move then the complexes $(\cc
**D,d)$ and $(\cc **{D^\prime},d^\prime)$ are homotopy equivalent.
\end{prop}

\begin{figure}[h]
\centerline{
\psfig{figure= 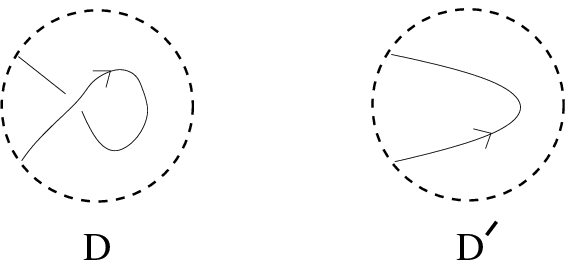}
}
\caption{}
\label{fig:ddprime} 
\end{figure}

 We are not going to provide a complete proof of this by any
 means. The aim is to give you an idea of how things go.  Let us look
 at the first Reidemeister move for a positive twist, so that diagrams
 $D$ and $D^\prime$ are identical except within a small region where
 they are shown in Figure \ref{fig:ddprime}.

\begin{figure}[h]
\centerline{
\psfig{figure= 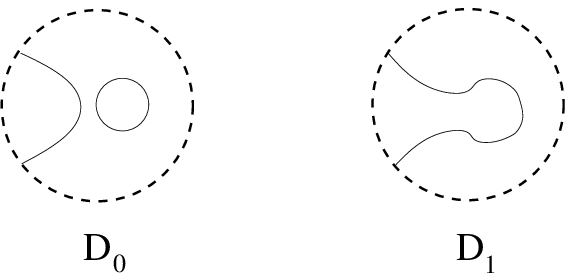}
}
\caption{}
\label{fig:d0d1} 
\end{figure}

We need to define chain maps $F\colon \cc ** {D^\prime} \ra \cc **D$
and $G\colon \cc **D \ra \cc ** {D^\prime}$ such that $GF-I$ and
$FG-I$ are null-homotopic. The first thing is to notice that 
we can split the vector space $\cc i*D$ as
\[
\cc i*D = \cc i*{D_0} \oplus \cc {i-1}*{D_1}
\]
where $D_0$ and $D_1$ are diagrams identical to $D$ except within the small regain where they are shown in Figure \ref{fig:d0d1}.

\begin{exercise}
Understand why there is the shift for $D_1$ in the above decomposition. What happens to the $q$-grading?
\end{exercise}

The differential $d$ can be written with respect to this splitting as
a matrix $\begin{pmatrix}d_0 & 0 \\ \delta & d_1 \end{pmatrix}$.

\begin{exercise}
Describe the map $\delta\colon \cc i*{D_0} \ra \cc {i-1+1}*{D_1}$ both algebraically and geometrically.
\end{exercise}

To define $F\colon \cc ** {D^\prime} \ra \cc **D$ we need to define
two coordinate maps $F_0 \colon \cc ** {D^\prime} \ra \cc **{D_0}$ and $F_1
\colon \cc ** {D^\prime} \ra \cc **{D_1}$ and then set $F=(F_0, F_1)$.
How should we go about constructing maps $\cc **{D^\prime} \ra \cc
**{D_0}$? Given a smoothing $\alpha^\prime$ of $D^\prime$ there is a
corresponding one $\alpha$ of $D_0$ (the one that resolves the
crossings in the same way). These smoothings look identical outside
the small region above. We can construct a cobordism from $\alpha^\prime$ to
$\alpha$ by taking a product with $I$ outside the small region and
inserting \inlinebigpic{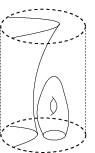} in the missing tube. By applying the TQFT to
this cobordism we get a map $V_{\alpha^\prime} \ra V_{\alpha}$. As
we can do this for each smoothing these maps assemble into a map $\cc
**{D^\prime} \ra \cc **{D_0}$. 

We can define another such map by gluing  \inlinebigpic{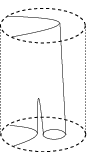} into the
missing tube. The map $F_0$ is the difference between the two maps
just defined. In pictures $F_0$ is the map defined by

\begin{figure}[h]
\centerline{
\psfig{figure= 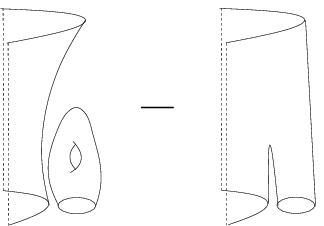}
}
\end{figure}

We take $F_1=0$ and set $F=(F_0,F_1)$.

As with Frobenius algebras we could write out $F$ more
algebraically if we wanted to. For each smoothing $\alpha$ of $D_0$
the vector space $V_{\alpha}$ is of the form $V_\alpha= Y_\alpha
\otimes V \otimes V$ - the last copy of $V$ for the separate circle we
see in the picture and the other copy of $V$ for the other circle
appearing. The corresponding smoothing of $D^\prime$ has associated to
it the vector space $Y_{\alpha^\prime} \otimes V$ - the copy of V for
the circle which enters the region shown. 
In this language $F_0\colon \cc **{D^\prime}\ra \cc **{D_0}$ is the map 
\[
F_0(y\ot v) = y\ot v \ot 2x - y \ot \Delta (v)
\]
(Remember here that $x$ is the degree -1 generator of $V$).

\begin{exercise}
Show that the bi-degree of $F$ is $(0,0)$.
\end{exercise}

\begin{figure}[h]
\centerline{
\psfig{figure= 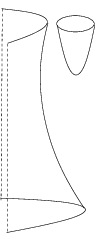}
}
\caption{}
\label{fig:tube4} 
\end{figure}

Now we turn to $G$ where we can be briefer. Define $G_0\colon \cc ** {D_0} \ra \cc **{D^\prime}$ by Figure \ref{fig:tube4} (using the method above) and let 
$G_1=0$..

\begin{exercise}
So far $F$ and $G$ are maps of vector spaces: check they are chain maps.
\end{exercise}

We now claim that $G$ and $F$ are part of a homotopy equivalence i.e. that $GF-I$ is null homotopic and $FG-I$ is null homotopic.

\begin{figure}[h]
\centerline{
\psfig{figure= 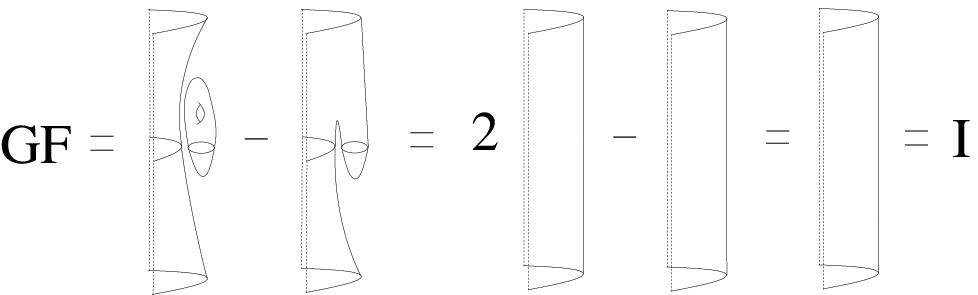}
}
\caption{}
\label{fig:tube5} 
\end{figure}

For the first of these we claim $GF=I$ (showing $GF-I$ is null homotopic via a trivial homotopy). This is where
using pictures comes into its own: the picture for the composition
$GF$ is simply gotten by placing one picture on top of the other as seen in Figure \ref{fig:tube5}.

Next we claim there is a map $H\colon\cc ** D \ra \cc{*-1}*D$ such
that $FG-I=Hd+dH$. Using the splitting above $H$ is the matrix
$\begin{pmatrix} 0 & h \\ 0 & 0\end{pmatrix}$ where $-h\colon \cc ** {D_1} \ra \cc ** {D_0}$ is
the map gotten by the method above using the picture

\begin{figure}[h]
\centerline{
\psfig{figure= 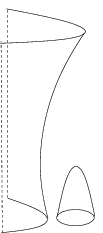}
}
\end{figure}

We compute $Hd + dH = \begin{pmatrix} h\delta & h d_1 + d_0h \\
 0 & \delta h\end{pmatrix}. $ Thus we need to show
\begin{eqnarray}
 h\delta  & = & F_0G_0 - I  \label{eq:HddH1}\\
 h d_1 + d_0h& =& 0 \label{eq:HddH2} \\
  \delta h & = & -I  \label{eq:HddH3}
\end{eqnarray}

\begin{figure}[h]
\centerline{
\psfig{figure= 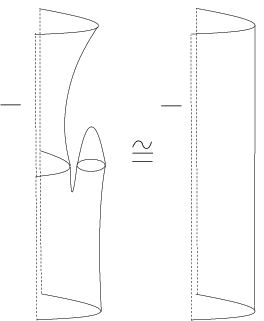}
}
\caption{}
\label{fig:tube7} 
\end{figure}

Equation (\ref{eq:HddH3}) is easy: just compose pictures as shown in Figure \ref{fig:tube7}.

Equation (\ref{eq:HddH2}) is essentially just due to the fact that
$d_0$ and $d_1$ are defined using the same cobordism (for $d_0$ there
is an extra cylinder) and by then looking carefully at the signs in
the definition of the differential.

\begin{exercise}
Check the assertions of the previous sentence.
\end{exercise}

Pictorially ({\ref{eq:HddH1}) is shown in Figure \ref{fig:tube8}. I know of no enlightening way to do it, but it is a
simple matter to see that this holds for $V$. (Remember the cobordism
is the identity outside the region so we need to check the
equality for maps $V\otimes V \ra V\ot V$.)

\begin{figure}[h]
\centerline{
\psfig{figure= 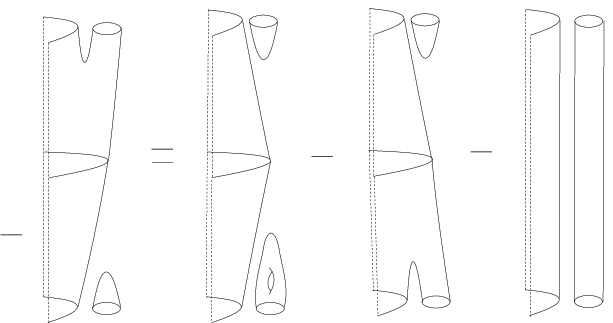}
}
\caption{}
\label{fig:tube8} 
\end{figure}

This concludes invariance under Reidemeister I positive twist: we have
produced $F$ and $G$ such that $FG-I$ and $GF-I$ are null-homotopic,
thus demonstrating that there is a homotopy equivalence $\cc **
{D^\prime} \simeq \cc **D$.

The above essentially follows Bar-Natan's proof - though Bar-Natan is
cleverer still: in his set-up one constructs a {\em geometric} complex
and works with tangles. One proves invariance without ever applying a
TQFT. This gives rise to a {\em universal theory} - more on this in
the next lecture. Refer to the end of the lecture for further remarks
and a reference.

\section{Khovanov homology}

Given an oriented link diagram $D$ we now define the {\em Khovanov homology
of the diagram} $D$ by
\[
\kh ** D = H(\cc ** D,d).
\]
By the previous section if $D$ is related to $D^\prime$ by a series of
Reidemeister moves then there is an isomorphism $\kh ** D \cong \kh **
{D^\prime}$. Thus if $L$ is an oriented link it makes sense to talk
about the {\em Khovanov homology of the link} $L$ (defined up to isomorphism
as the Khovanov homology of any diagram representing it).

\begin{prop}
\[
\sum (-1)^i \qdim {\kh i*L} = \hat J (L)
\]
\end{prop}

\begin{proof}
It is an exercise in linear algebra to show that \[
\sum (-1)^i \qdim
{\kh i*D} = \sum (-1)^i \qdim {\cc i*D} 
\]
 and we have already observed
that the right-hand side  is $\hat J (L)$.
\end{proof}

Khovanov homology is a stronger invariant than the Jones polynomial as
the following example illustrates.
\begin{exe}\label{ex:khstronger}
\end{exe}

{\em
\begin{center}
\begin{tabular}{lll}
$D_1 =$ \raisebox{-7mm}{\psfig{figure={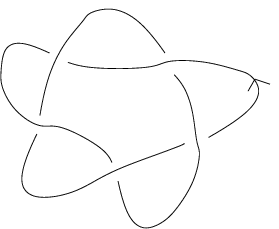}}} $\;\;\;\;$ & $\;\;\;\;\;\;\;\;\;\;\;\;\;\;\;$ &$D_2=$\raisebox{-7mm}{\psfig{figure={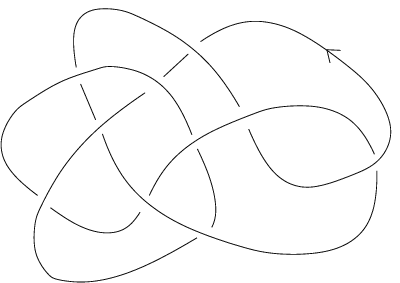}}} \\
\end{tabular}
\end{center}
}

{\em
\begin{center}
\begin{tabular}{lll}
Unnormalised Jones polynomial: && Unnormalised Jones polynomial\\
 $\hat J (D_1) = q^{-3} + q^{-5} + q^{-7} - q^{-15}$ & & $\hat J (D_2)=  q^{-3} + q^{-5} + q^{-7} - q^{-15}$\\
 && \\
Khovanov homology $\kh ij {D_1}$: && Khovanov homology $\kh ij {D_2}$:\\
&&\\
\begin{tabular}{|c||c|c|c|c|c|c|}
\hline
\backslashbox{$j$}{$i$} &-5&-4&-3& -2 & -1& 0\\
\hline\hline
-3 & &&&&& $\bQ$\\\hline
-5 &&&&&& $\bQ$\\\hline
-7 &&&& $\bQ$&&\\\hline
-9&&&&&&\\\hline
-11&&$\bQ$&$\bQ$&&&\\\hline
-13&&&&&&\\\hline
-15&$\bQ$&&&&&\\\hline
\end{tabular}
&&
\begin{tabular}{|c||c|c|c|c|c|c|c|c|}
\hline
\backslashbox{$j$}{$i$}&-7&-6&-5&-4&-3& -2 & -1& 0\\
\hline\hline
-1&&&&&&& $\bQ$& $\bQ$\\\hline
-3 & &&&&&&& $\bQ$\\\hline
-5 &&&&& $\bQ$ &  $\bQ\oplus\bQ$&&\\\hline
-7 &&&& $\bQ$&&&&\\\hline
-9&&&& $\bQ$& $\bQ$&&&\\\hline
-11&&$\bQ$& $\bQ$&&&&&\\\hline
-13&&&&&&&&\\\hline
-15&$\bQ$&&&&&&&\\\hline
\end{tabular}
\end{tabular}
\end{center}
}

\vspace*{3mm}

The missing rows in the above table (those with even $q$-degree) are
all trivial. In fact this is a more general phenomenon.

\begin{prop}
If a link $L$ has an odd number of components then $\kh *
{\text{even}} L =0$. If $L$ has an even number of components then $\kh
*{\text{odd}} L =0$.
\end{prop}

\section{Notes and further reading} 
\noindent
To find out more about Frobenius algebras and TQFTs two good places to start are \cite{kock} and \cite{abrams}. For a more general treatment of TQFTs consult
\cite{turaev}.

\vspace*{2mm}

\noindent
Reidemeister invariance was first proved by Khovanov in his original paper
\cite{khovanov1} and proofs can also be found in \cite{barnatan2} which also 
contains computations of the Khovanov homology of prime knots with
diagrams with up to ten crossings. This is where the computation in
Example \ref{ex:khstronger} is taken from. There are some simple examples in
\cite{wehrli} showing that there exist mutant links (and so having the 
Jones polynomial) which are separated by Khovanov homology. At the time of writing 
it is unknown whether mutant {\em knots} can be separated by Khovanov homology.

\vspace*{2mm}

\noindent
The proof of invariance under Reidemeister move I presented above is
closer to the proof found in \cite{barnatan1}. In this paper Bar-Natan
considers the cube (with smoothings at vertices and cobordisms on
edges) as a {\em geometric complex} (i.e. a complex in an abelianized
category of cobordisms). In order to prove invariance (in the homotopy
category of these geometric complexes) one needs to take a quotient of
the cobordism category by the relations shown in Figure
\ref{fig:barnatanrels}.
\begin{figure}[h]
\centerline{
\psfig{figure= 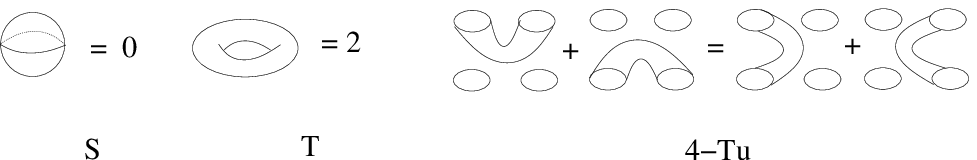}
}
\caption{}
\label{fig:barnatanrels} 
\end{figure}
It is the {\em 4-Tu} relation which is used in Figure
\ref{fig:tube8}. Things are better even than this: the theory is
completely local and one works with tangles (we had pictures of
tangles, but kept in our minds the fact that these were part of a
larger diagram).



\chapter{Lecture Three}
In this lecture we begin by looking at a long exact sequence in
Khovanov homology. Then we examine the kind of functoriality present
and briefly discuss the invariants of embedded surfaces in $\bR^4$
thus defined. We end with a look at theories defined over different
base rings.

\section{A long exact sequence}
In algebraic topology there are many theoretical tools for computation
such as long exact sequences, spectral sequence and so on. In Khovanov
homology there is less available in the arsenal, but there is one
useful long exact sequence which we now discuss.

If we choose a crossing of a diagram $D$ we can  resolve it in the two possible ways to give two new diagrams $D_0$ and $D_1$ as in Figure \ref{fig:twoways}.

\begin{figure}[h]
\centerline{
\psfig{figure= 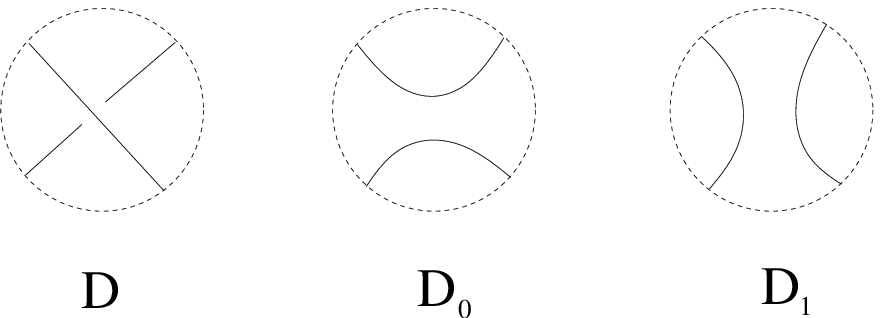}
}
\caption{}
\label{fig:twoways} 
\end{figure}

Ignoring gradings for a moment there is a decomposition of (ungraded) vector spaces
\[
C(D) = C(D_0) \oplus C(D_1).
\]

In fact $C(D_1)$ is a sub-complex and there is a short exact sequence
\[
0\ra C(D_1) \ra C(D) \ra C(D_0) \ra 0.
\]
Putting the gradings back in requires a bit of care though it is not hard.

{\bf Case I:} the selected crossing of $D$ is a negative crossing. In
this case $D_1$ inherits an orientation from $D$ (the 1-smoothing is
the orientation preserving smoothing). For $D_0$ there is no orientation consistent with $D$ so just orient it as you please. Let 
\[
c = \text{ number of negative crossings in $D_0$} -  \text{ number of negative crossings in $D$}.
\]
Then for each $j$ there is a short exact sequence
\[
0\ra \cc i {j+1} {D_1} \ra \cc ijD \ra \cc {i-c}{j-3c-1}{D_0} \ra 0
\]
and hence a long exact sequence
\[
\xymatrix{
\ar[r]^-{\delta_*} &\kh i {j+1} {D_1} \ar[r] & \kh ijD \ar[r]&  \kh {i-c}{j-3c-1}{D_0} \ar[r]^-{\delta_*} & \kh {i+1} {j+1} {D_1}\ar[r]&  .
}
\]

If we write the differential of $\cc ** D$ as a matrix
$\begin{pmatrix}d_0 & 0 \\ \delta & d_1 \end{pmatrix}$ then the
boundary map in the long exact sequence is $\delta_*$ i.e. the map
induced in homology (suitably shifted to take account of the new
orientation of $D_0$).

{\bf Case II:} the selected crossing of $D$ is a positive crossing. In
this case $D_0$ inherits an orientation from $D$ (this time the
0-smoothing is the orientation preserving smoothing). For $D_1$ there
is no orientation consistent with $D$ so just orient it as you
please. Let
\[
c = \text{ number of negative crossings in $D_1$} -  \text{ number of negative crossings in $D$}.
\]
Then for each $j$ there is a short exact sequence
\[
0\ra \cc {i-c-1} {j-3c-2} {D_1} \ra \cc ijD \ra \cc {i}{j-1}{D_0} \ra 0
\]
and hence a long exact sequence
\[
\xymatrix{
\ar[r]^-{\delta_*} &\kh {i-c-1} {j-3c-2} {D_1} \ar[r] & \kh ijD \ar[r]&  \kh {i}{j-1}{D_0} \ar[r]^-{\delta_*} & \kh {i-c} {j-3c-2} {D_1}\ar[r]&.
}
\]

\begin{exercise}
Check the gradings in the long exact sequences above.
\end{exercise}

\begin{exe}
At the end of Lecture 1 we computed the Khovanov homology of the Hopf
link. Let us re-do this calculation using the Khovanov homology of the
unknot and the long exact sequence. We choose to resolve the top
crossing (a negative crossing) thus giving $D_0$ and $D_1$ as shown in
Figure \ref{fig:hopfresolve}

\begin{figure}[h]
\centerline{
\psfig{figure= 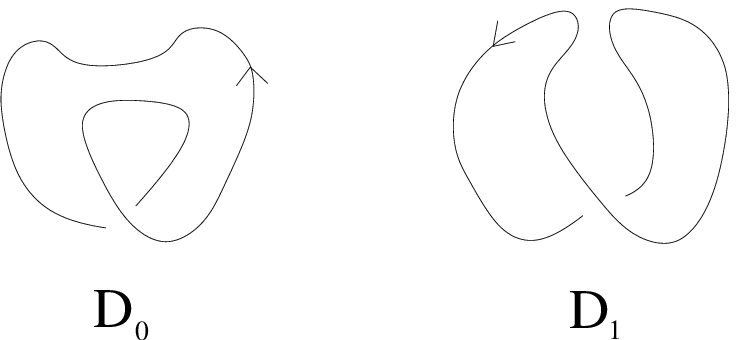}
}
\caption{}
\label{fig:hopfresolve} 
\end{figure}

Here we have $c=-2$ so the long exact sequence is
\[
\xymatrix{
\ar[r]^-{\delta_*} &\kh i {j+1} {D_1} \ar[r] & \kh ijD \ar[r]&  \kh {i+2}{j+5}{D_0} \ar[r]^-{\delta_*} & \kh {i+1} {j+1} {D_1} \ar[r] &.
}
\]
Since $D_0$ and $D_1$ are both the unknot they only have non-trivial homology in homological degree 0 (where there are generators in $q$-degree $+1$ and $-1$). Thus the long exact sequence breaks up and there are two interesting pieces
\begin{equation}\label{eq:hlles1}
0\ra \kh 0{j+1}{D_1} \ra \kh 0jD \ra 0 \ra 0
\end{equation}
\begin{equation}\label{eq:hlles2}
0\ra 0\ra \kh {-2}{j}{D} \ra \kh 0{j+5}{D_0} \ra 0 
\end{equation}
From (\ref{eq:hlles1}) we see that all groups are zero unless $j=0,-2$
from which we conclude $\kh 00D \cong \kh 0{-2}D \cong
\bQ$. Similarly, from (\ref{eq:hlles2}) we see that all groups are
zero unless $j=-4,-6$ from which we conclude $\kh {-2}{-4}D \cong \kh
{-2}{-6}D \cong \bQ$. This result is happily in agreement with the
computation at the end of Lecture 1.
\end{exe}

\section{Functorial properties}

It is convenient to study links by projecting onto the plane and
studying link diagrams instead. The diagrammatic representation of a
given link is far from unique, but this is well understood: two
diagrams represent the isotopic links if and only if they are related by
Reidemeister moves.

Something similar is true for link cobordisms (recall that a link
cobordism $(\Sigma, L_0,L_1)$ is a smooth, compact, oriented surface
$\Sigma$ generically embedded in $\bR^3\times I$ such that
$\partial\Sigma = \ol{L}_0\sqcup L_1$ with $\partial\Sigma \subset
\bR^3\times \{0,1\}$). A link cobordism can be represented by a
sequence of oriented link diagrams - the first in the sequence being a
diagram $D_0$ for $L_0$ and the last being a diagram $D_1$ for
$L_1$. Two consecutive diagrams in this sequence must be related by a
small set of allowable moves which are
\begin{enumerate}
\item Reidemeister I, II or III moves,
\item Morse 0-,1- or 2-handle moves.
\end{enumerate}
The Reidemeister moves are just the usual ones and the Morse moves are shown in Figure \ref{fig:morsemoves}.

\begin{figure}[h]
\centerline{
\psfig{figure= 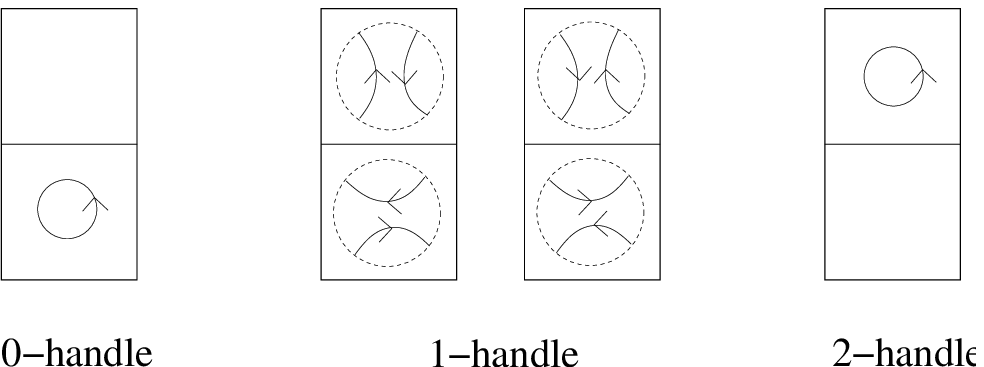}
}
\caption{}
\label{fig:morsemoves} 
\end{figure}

Geometrically the Morse moves are shown in Figure \ref{fig:morsemovesgeom}. 

\begin{figure}[h]
\centerline{
\psfig{figure= 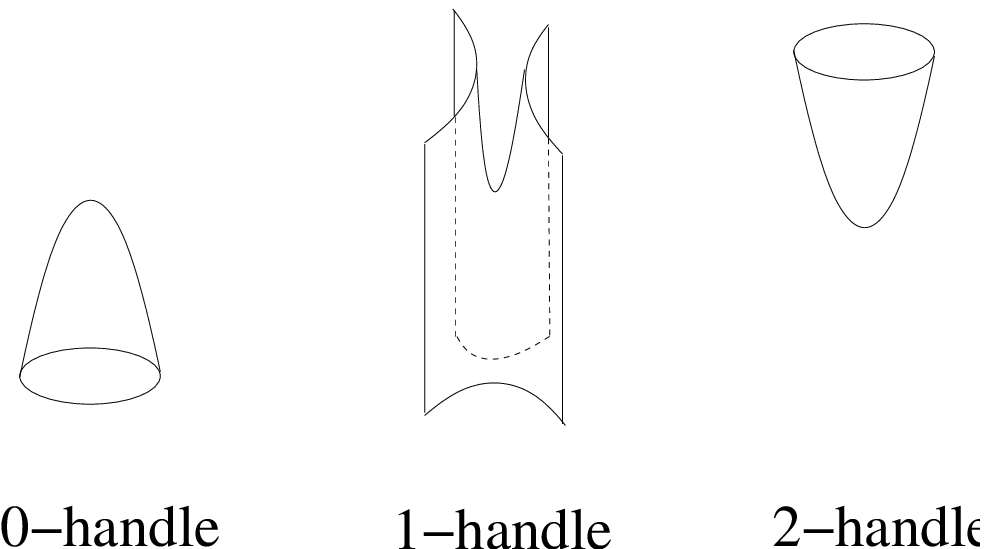}
}
\caption{}
\label{fig:morsemovesgeom} 
\end{figure}

Such a sequence of diagrams is known as a {\em movie}. 

\begin{exe}
Figure \ref{fig:movie} is a movie representing a cobordism from the
Hopf link to the empty cobordism (drawn across the page rather than
down to save space).
\end{exe}

\begin{figure}[h]
\centerline{
\psfig{figure= 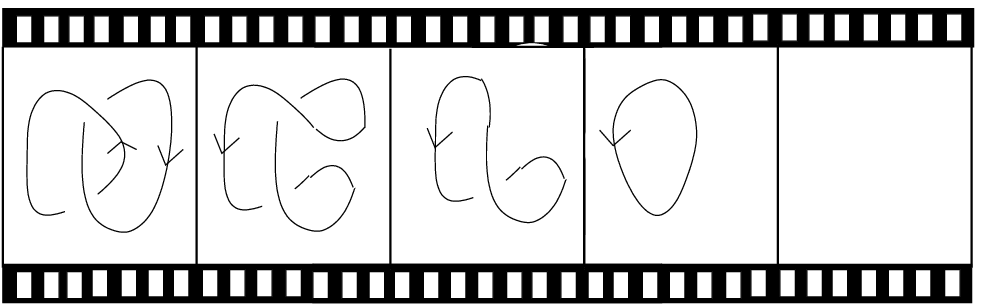}
}
\caption{}
\label{fig:movie} 
\end{figure}

A movie representation of a link cobordism is not unique: there may be
many different movies of the same cobordism. However, again this is
well understood: two movies represent isotopic link cobordisms if and
only if they are related by a series of {\em movie moves} or by
interchanging the levels of distant critical points. Each movie move
replaces a small clip of the movie by a different clip.  We will not
go into this in more detail here.

A movie $(M,D_0,D_1)$ induces a map on Khovanov homology
\[
\phi_M\colon \kh ** {D_0} \ra \kh *{*+\chi}{D_1}
\]
in the following way. (Here $\chi$ is  the number of Morse 0- and 2-handle moves minus the number of Morse 1-handle moves).
We will define a map between each two consecutive frames of the movie
and then compose all of these to get $\phi_M$. For Reidemeister moves
we have already argued that there is a homotopy equivalence of chain
complexes which gives a map in homology and it is this map we
take. For a 0-handle move, if the before-frame consists of a link diagram
$D$ then the after frame consists of $D \sqcup
\text{unknot}$. Since $\kh ** {D \sqcup \text{unknot}} = \kh **D
\otimes V$ we take the the map $\kh ** D \ra \kh **D \otimes V$ to be
$Id \otimes i$ where $i\colon \bQ \ra V$ is the unit of the Frobenius
algebra $V$. Since $1\in V$ has $q$-degree 1 this map increments $q$-degree by one. For the 2-handle move we do a similar thing using the
counit of the Frobenius algebra. 

For the 1-handle move let $D$ and $D^\prime$ be the before- and
after-frames of the move.  We construct a map $\cc **D \ra \cc
*{*-1}{D^\prime}$ by using the geometric techniques at the beginning
of Lecture 2. For each smoothing $\alpha$ of $D$ there is a
corresponding one $\alpha^\prime$ of $D^\prime$ different only in the
small region in which the move takes place. A cobordism can be
constructed from $\alpha$ to $\alpha^\prime$ by taking a product with
$I$ outside the small region and inserting a saddle in the missing
tube. Do this for each smoothing, apply the TQFT, assemble the
resulting maps and take homology to get a map $\kh **D \ra \kh
*{*-1}{D^\prime}$.

\begin{prop}
If $(M,D_0,D_1)$ is related to $(M^\prime, D_0,D_1)$ by a sequence of
movie moves or interchanging the levels of distant critical points
then $\phi_{M^\prime} = \pm \phi_M$.
\end{prop}

We will not prove this theorem. The sign discrepancy is annoying but
some movie moves (though not all!) change the sign. In order to say
that $\kh ** -$ is a functor you therefore need to projectivize the
target category.

\section{Numerical invariants of closed surfaces}
Using the above one can define a numerical invariant of closed
oriented surfaces smoothly embedded in $\bR^4$. Such a surface $\Sigma$ may be
regarded as a link cobordism between the empty link and the empty
link. Thus, representing $\Sigma$ by a movie $M$ and noting that $\kh
** \emptyset = \bQ$, the above discussion gives us a map $\phi_M\colon
\bQ \ra \bQ$. The {\em Khovanov-Jacobsson number}, of the embedded surface $\Sigma$ is defined to be $KJ_\Sigma = |\phi_M(1)|$.

Unfortunately, these numbers are rather disappointing. If
$\chi(\Sigma)$ is non-zero then $KJ_\Sigma = 0$ since $\phi_M$ shifts
the $q$-degree by $\chi(\Sigma)$ (and both the source and target of $\phi_M$ are non-zero only in bi-degree $(0,0)$). For embedded tori there is the
following result.

\begin{prop}
If $\Sigma$ is a smoothly embedded torus in $\bR^4$ then $KJ_\Sigma=2$.
\end{prop}

\section{Coefficients and Torsion}
So far we have been working over the rational numbers, but nothing we
have said so far really relies on this. Everything remains valid (the
construction of a complex, the proofs of invariance etc) replacing
$\bQ$ by any commutative ring with unit. Instead of ``vector space''
you need to write ``projective $R$-module of finite type''. In
particular you can work over the integers and ask if, like for
ordinary homology of spaces, the interesting phenomenon of {\em
torsion} emerges. It does.

The first place this is seen is for the trefoil
\inlinepic{trefoil.eps}. The cube of this trefoil (as requested in the
exercise at the end of Lecture 1) is given in Figure
\ref{fig:trefoilcube}.

\begin{figure}[h]
\centerline{
\psfig{figure= 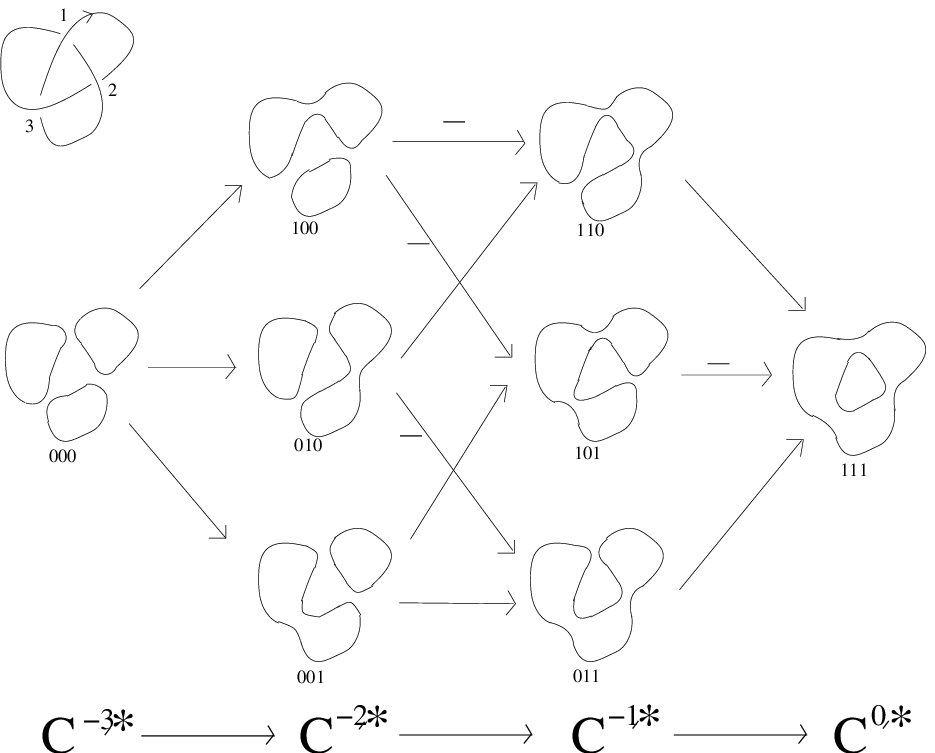}
}
\caption{}
\label{fig:trefoilcube} 
\end{figure}

In bi-degree $(-2,-7)$ the cycles are generated by 
\[
\{z_1=(x\ot
x,0,0),z_2=(0,x\ot x,0), z_3=(0,0,x\ot x)\}
\]
 and in bi-degree
$(-3,-7)$ the chains are generated by 
\[
\{ c_1=(1,x,x),c_2=(x,1,x),
c_3=(x,x,1)\}.
\]
Recalling the definition of the differential one easily sees
\[
d(c_1) = z_1 + z_3 \;\;\;\; 
d(c_2) = z_2 + z_3 \;\;\;\; 
d(c_3) = z_1 + z_2
\]
Thus in homology we have $[z_1]=[z_2]=[z_3]$. Note also that
\[
d(c_1+c_3-c_2) = z_1+z_3+z_1+z_2-z_2-z_3 = 2z_1.
\]
Thus, rationally $z_1$ is a boundary (it is hit by
$\frac{1}{2}(c_1+c_3-c_2)$) and our potential homology class above is
trivial. Over the integers $[z_1]$ is a non-trivial homology class,
but $2[z_1]$ is trivial, thus in homology we have a copy of $\bZ/2$. 

For the record the full integral homology of the trefoil
\inlinepic{trefoil.eps} is given below.
\begin{center}
\begin{tabular}{|c||c|c|c|c|}
\hline
\backslashbox{$j$}{$i$} &-3& -2 & -1& 0\\
\hline\hline
-1 &&&& $\bZ$\\\hline
-3 &&&& $\bZ$\\\hline
-5 &&$\bZ$&& \\\hline
-7 && $\bZ/2$&&\\\hline
-9&$\bZ$ &&&\\\hline
\end{tabular}
\end{center}

In fact torsion abounds as shown in the following result (which we do not prove):

\begin{prop}\label{prop:torsion}
The integral Khovanov homology of every alternating link, except the
trivial knot, the Hopf link and their connected sums and disjoint
unions, has torsion of order two.
\end{prop}

In order to describe Khovanov homology with coefficients in a ring $R$
in terms of integral Khovanov homology we apply a standard result in
homological algebra: the universal coefficient theorem. Since $\ccc
**DR = \ccc **D {\bZ} \otimes_\bZ R$ the universal coefficient theorem
tells us that there is a short exact sequence
\[
\xymatrix{
0 \ar[r] & \khc ijD{\bZ} \otimes_{\bZ} R \ar[r] & \khc ijD {R} \ar[r] & \Tor (\khc {i+1}jD{\bZ} , R)\ar[r]& 0 .
}
\]

\begin{exercise}
Use the short exact sequence above to compute $\khc **{\inlinepic{trefoil.eps}} {\bZ/2}$.
\end{exercise}

Closely related to this is the K\"unneth formula, which we can use to compute the Khovanov homology of a disjoint union. Given two link diagrams $D_1$ and $D_2$ then $\cc **{D_1\sqcup D_2} \cong \cc** {D_1} \otimes \cc **{D_2}$ or more precisely:
\[
\cc i j {D_1\sqcup D_2} \cong \bigoplus_{\substack{p+q=i\\s+t=j}}\cc ps {D_1} \otimes \cc qt{D_2}.
\]
Thus the K\"unneth formula gives us a split short exact sequence
\begin{align*}
0 \ra \bigoplus_{\substack{p+q=i\\s+t=j}}\khc ps{D_1}R \otimes & \khc qt{D_2}R
\ra  \khc ij {D_1\sqcup D_2}R \\ & \ra  \bigoplus_{\substack{p+q=i+1\\s+t=j}}\Tor^R_1 (\khc ps{D_1}R, \otimes \khc qt{D_2}R)\ra 0.
\end{align*}
Over $\bQ$ the Tor group is always trivial so we have
\[
 \khc ij {D_1\sqcup D_2}{\bQ} \cong  \bigoplus_{\substack{p+q=i\\s+t=j}}\khc ps{D_1}{\bQ} \otimes \khc qt{D_2}{\bQ}.
\]

 \section{Notes and further reading} 
\noindent
The long exact sequence is implicit in Khovanov's original paper, but appeared in a slightly different form in \cite{viro}. Lee used the singly graded version in \cite{lee}.
It has appeared in a variety of places since then and with gradings as we have given them in \cite{rasmussen4}. 
This is an interesting survey paper in its own right discussing parallels between Khovanov homology and knot Floer homologies.

\vspace*{2mm}

\noindent
You can read about, cobordisms and their representations as movies in \cite{cartersaito}.

\vspace*{2mm}

\noindent
Khovanov conjectured the functoriality properties in his original paper \cite{khovanov1}. This was then proved in \cite{jacobsson}
and independently in \cite{khovanov2}.
Using his geometric techniques Bar-Natan proved functoriality (in more generality) in \cite{barnatan1}.

\vspace*{2mm}

\noindent
The proposition about Khovanov-Jacobsson numbers has been proved for a certain class of torus embeddings in \cite{cartersaitosatoh}
and then for all torus embeddings in \cite{tanaka} and independently using different techniques in \cite{rasmussen2}.

\vspace*{2mm}

\noindent
Torsion in Khovanov homology has been studied in a number of
papers. The best places to start would be \cite{shumakovitch} and
\cite{asaedaprzytycki}.  Both of these prove Proposition
\ref{prop:torsion} concerning 2-torsion and the former has a number of
interesting conjectures about torsion. One of these conjectures that
all torsion is 2-torsion, which is now known to be false. Bar-Natan's
computer program calculates $\kh {22}{73} {T(8,7)} = \bZ/2 \oplus
\bZ/4 \oplus \bZ/5 \oplus \bZ/7$, where $T(8,7)$ is the torus link with 7 strands and 8 positive twists.



\chapter{Lecture Four}
In this lecture we begin by discussing a family of Khovanov-type link
homology theories, then focus on the first of these to be defined, Lee
theory. This naturally leads to the work of Rasmussen on a new
concordance invariant of knots which has many wonderful properties.

\section{A family of Khovanov-type theories}
An obvious question to ask is: can we replace the Frobenius algebra
$V$ used to construct the Khovanov complex by some other Frobenius
algebra and still get a link homology theory (i.e. an invariant with
nice functorial properties)? Or, re-phrased, what conditions must a
Frobenius algebra $A$ satisfy to give a link homology theory?

For simplicity let us again take $\bQ$ as the base ring. It is relatively easy to see that we must have dim$(A)=2$. Consider the two representations of the unknot below.

\begin{figure}[h]
\centerline{
\psfig{figure= 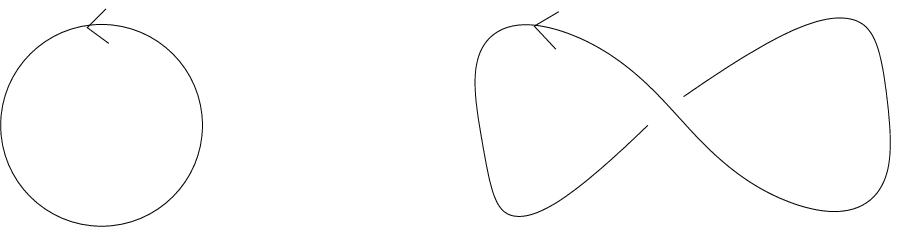}
}
\end{figure}

The first gives a complex
\[
0\ra A \ra 0 
\] 
and the second a complex
\[
0\ra A \ra A\ot A \ra 0.
\]
In the first, the copy of $A$ is in degree zero and in the second in degree -1. Since both diagrams represent the unknot we require these two  complexes to be homotopy equivalent. The Euler characteristic of homotopy equivalent complexes must be equal, hence
\[
\text{dim}(A) = - \text{dim}(A) + \text{dim}(A\otimes A)
\]
so $2 \text{dim}(A) - \text{dim}(A)^2 = 0$ and we conclude $\text{dim}(A)=2$.

It turns out that additionally one needs $\epsilon(1) =0$ and this is
all that is required for a Frobenius algebra to give rise to a link
homology theory. To understand why this is the case one needs to use
Bar-Natan's geometric theory. As explained in the ``Further reading''
of Lecture 2, Bar-Natan imposes three relations on the category of
cobordisms: {\em S,T} and {\em 4-Tu}. If we wish to apply a TQFT to
his geometric setting, the underlying Frobenius algebra $A$ must
satisfy these relations (or at least their algebraic
counterparts). The relation {\em S} just says $\epsilon(1)=0$ and {\em
T} is the condition dim$(A)=2$. While the {\em 4-Tu} relation is
necessary geometrically it is automatically satisfied for two
dimensional Frobenius algebras satisfying $\epsilon(1)=0$. This is an
indication of the power of Bar-Natan's geometric approach - things are
proved in terms of {\em complexes of cobordisms} and so applying a
different TQFT (satisfying the necessary conditions) simply gives
another theory with no additional effort.

As a vector space we may as write $A=\bQ\{1,x\}$ as before. Then
normalising so that $\epsilon(x)=1$ we have a family of theories one
for each pair $(h,t)\in \bQ\times \bQ$. The Frobenius algebra
$A_{h,t}$ has multiplication given by
\[
1^2=1 \spaces 1x=x1= x \spaces x^2=hx+t1
\]
and comultiplication given by
\[
\Delta(1) = 1\ot x + x\ot 1 -h 1\ot 1 \spaces \Delta(x) = x\ot x + t 1\ot 1.
\]
The unit and counit are
\[
i(1) = 1 \spaces \epsilon(1) =0\spaces \epsilon(x) =1.
\]
A link homology theory is obtained by carrying out exactly the same
construction as outlined in Lecture 1, but replacing the Frobenius
algebra $V$ with $A_{h,t}$. When $h=t=0$ then you get the original
theory of Lecture 1.

The first such variant to be studied was the case $(h,t) = (0,1)$ by
E.S. Lee giving a theory now known as {\em Lee Theory} - the topic of
the next section.

In fact, up to isomorphism Lee theory is the only other rational theory in the family under discussion. 

\begin{prop}\label{prop:classoverq}
If $h^2+4t =0$ then the resulting theory is isomorphic to the original
Khovanov homology and if $h^2+4t\neq 0$ then the resulting theory is
isomorphic to Lee theory.
\end{prop}

\section{Lee theory}
The alert reader will have noticed a minor problem in carrying out the
construction in Lecture 1 using the Frobenius algebra
$A_{h,t}$. This is that one loses the $q$-grading. The degree of
$x^2=hx+t1$ is not even homogeneous if $h$ and $t$ are both
non-zero. In fact the only case where a second grading exists is
$h=t=0$. For now let us simply ignore the $q$-grading: all theories
will be singly graded by the homological grading. In the next section
we will see that in fact that we do not need to completely throw away
the second grading - it just gets replaced by a filtration instead.

Lee theory ($h=0, t=1$) was the first variant of Khovanov homology to
appear and remarkably it can be computed explicitly. 

\begin{prop}
The dimension of $\lee * L$ is $2^k$ where $k$ is the number of
components in $L$. 
\end{prop}

Things are even better still and there are explicit generators whose
construction is as follows. There are $2^k$ possible orientations of
$L$. Given an orientation $\theta$ there is a {\em canonical
smoothing} obtained by smoothing all positive crossings to
0-smoothings and all negative crossings to 1-smoothings. For this
smoothing one can divide the circles into two disjoint groups, Group 0
and Group 1 as follows. A circle belongs to Group 0 (Group 1) if it
has the counter-clockwise orientation and is separated from infinity
by an even (odd) number of circles or if it has the clockwise
orientation and is separated from infinity by an odd (even) number of
circles. Figure \ref{fig:groups} shows an orientation of the Borromean rings, its canonical smoothing and division into groups.

\begin{figure}[h]
\centerline{
\psfig{figure= 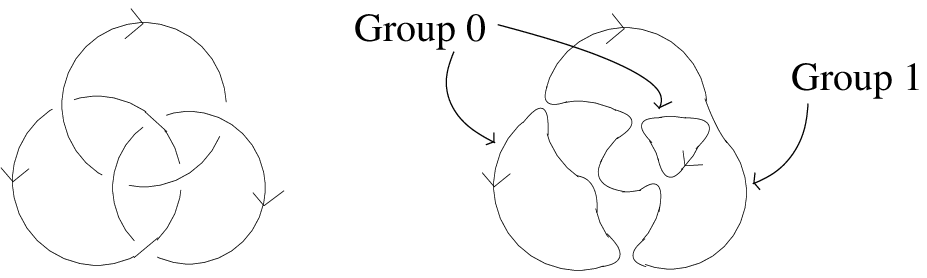}
}
\caption{}
\label{fig:groups} 
\end{figure}

Now consider the element in the chain complex for $L$ defined by
labelling each circle from Group 0 with $x+1$ and each circle from
Group 1 with $x-1$. It turns out that this defines a cycle, $\cano$,
and the homology class thus defined is a generator. Moreover all
generators are obtained this way and one has:
\[
\lee * L  \cong \bQ \{ [\cano] \mid \theta \mbox{ is an orientation of $L$}\}
\]  
This is not supposed to be obvious - read Lee's paper to find out why.

It is also possible to determine the degree of the generators in terms
of linking numbers. Let $L_1, \ldots , L_k$ denote the components of
$L$. Recalling that $L$ is oriented from the start, if we are given
another orientation of $L$, say $\theta$, then we can obtain $\theta$
by starting with the original orientation and then reversing the
orientation of a number of strands. Suppose that for the orientation
$\theta$ the subset $E\subset \{1, 2, \cdots , k \}$ indexes this set
of strands to be reversed. Let $\ol{E}=\{1,\ldots,k\}\backslash
E$. The degree of the corresponding generator $[\cano]$ is then given
by
\[
\mbox{deg}([\cano]) = 2 \times \sum_{l\in E,m\in \ol{E}} \lk(L_l,L_m)
\]
where $\displaystyle{\lk(L_l,L_m)}$ is the linking number (for the
original orientation) between component $L_l$ and $L_m$.

\begin{exercise}
Compute $\lee * {\inlinepic{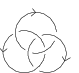}}$.
\end{exercise}

Since Lee theory is a link homology theory (and so has nice functorial properties) one can ask how canonical generators behave under cobordisms.

\begin{prop} \label{prop:cobgens}
Let $(\Sigma, L_0, L_1)$ be a cobordism presented by a movie
$(M,D_0,D_1)$. Suppose that every component of $\Sigma$ has a boundary
component in $L_0$. Then the induced map $\phi_M\colon \lee * {D_0}
\ra \lee *{D_1}$ has the property that $\phi_M([\canoo])$ is a non-zero multiple of $[\canooo]$, where
$\canoo$ and $\canooo$ are the orientations induced by the orientation
of $\Sigma$.
\end{prop}

\section{Rasmussen's invariant of knots}
Even though Lee theory is only singly graded it possesses a filtration
which can be used to define a new concordance invariant of
knots. Recall that originally we defined the $q$-grading of a chain
$v\in C^i(D)$ by
\[
q(v) = \text{deg}(v) + i + n_+ - n_-.
\]
In Lee theory we end up with elements that are not homogeneous with
respect to $q$-degree. However, for any monomial $w$ the quantity
$q(w)$ still makes sense and for an arbitrary element $v\in C^i(D)$
which can be written as a sum of monomials $v=v_1 + \cdots v_l$ we set
\[
q(v) = \text{min}\{q(v_i) \mid i = 1,\ldots , l\}.
\]
This defines a decreasing filtration on $C^*(D)$ by setting
\[
F^kC^*(D) = \{v\in C^*(D) \mid q(v) \geq k\}.
\]
The differential in $C^*(D)$ is a filtered map and thus Lee theory is a {\em filtered} theory.

Passing to homology we define for $\alpha \in \lee * D$
\[
s(\alpha) = \text{max}\{ q(v) \mid [v]=\alpha \}
\]
i.e. look at all representative cycles of $\alpha$ and take their
maximum $q$-value. Now for a knot $K$ define
\[
\smin = \text{min}\{s(\alpha) \mid \alpha \in \lee 0 K, \alpha \neq 0\},
\]
\[
\smax = \text{max}\{s(\alpha) \mid \alpha \in \lee 0 K, \alpha \neq 0\}
\]
and finally {\em Rasmussen's $s$-invariant of $K$} is 
\[
s(K)= \frac{\smin + \smax}{2}. 
\]

It turns out that $\smax = \smin +2$ and so $s(K)$ is always an
integer. We have the following properties.
\begin{enumerate}
\item $s(K)$ is an invariant of the concordance class of $K$,
\item $s(K_1\# K_2) = s(K_1) + s(K_2)$,
\item $s(K^!) = -s(K)$, where $K^!$ is the mirror image of $K$.
\end{enumerate}

These are not obvious, but we refer the reader to the original
reference for a proof.

In general it is hard to calculate the $s$-invariant of a knot. For
{\em positive knots} (one which has a diagram with only positive crossings)
it is easy. 
\begin{exe}
Let $K$ be a positive knot and $D$ a diagram for $K$. Since all the crossings are positive, there is only one
smoothing making up homological degree zero: the canonical
smoothing. Thus the canonical generator from the given orientation
lies in degree zero. Since $C^{-1}(D)=0$, the only representative of
$[\cano]$ is $\cano$ itself so $s([\cano]) = q(\cano)$.

The minimum possible $q$-value in degree zero is when each circle of
the canonical smoothing is labelled with $x$, and this does occur: as
a monomial in $\cano$. Thus $\smin = s([\cano]) = q(\cano) = -r +n$
where $r$ is the number of circles in the canonical smoothing. Thus
$s(K) = -r +n +1$.
\end{exe}

\begin{exercise}
Show that the $s$-invariant of the $(p,r)$-torus knot is $(p-1)(r-1)$.
\end{exercise}

One of the most interesting properties of $s$ is that it provides a
lower bound for the {\em slice genus} (also known as the {\em 4-ball
genus}). Recall that the slice genus $g^*(K)$ is the minimum possible genus of
a smooth surface-with-boundary smoothly embedded in $B^4$ with $K\subset \partial B^4$
as its boundary.

\begin{prop}
$|s(K)| \leq 2g^*(K)$
\end{prop}

We will prove this as it is relatively easy and is a great
demonstration of the usefulness of functoriality of link homology. Let
$\Sigma$ be a smooth surface of genus $g$ smoothly embedded in $B^4$
with boundary the knot $K$. We can remove a small disc from $\Sigma$
to get a smooth cobordism from $K$ to the unknot $U$. We can represent
this cobordism by a movie $(M, D, U)$ (here $D$ is a diagram for
$K$). Since $\lee * -$ has the functorial property described in Lecture
3 there is a map
\[
\phi_M\colon \lee * D \ra \lee * U = \bQ\{1,x\}.
\] 
It turns out this map has filtered degree
$\chi(\Sigma) = -2g$. 

Now let $\alpha\in \lee 0 K$ be a non-zero element such that $s(\alpha) = \smax$. Again by applying Proposition \ref{prop:cobgens} we get $\phi_M(\alpha)$ is non-zero in $\lee 0 U$ so
\[ 
1 = \smaxu \geq s(\phi_M(\alpha)) \geq s(\alpha) -2g = \smax -2g.
\]
Thus since $\smax=s(K) +1$ we have $s(K)\leq 2g$ and since this argument applies to any surface (including one of minimal genus) we get
\[
s(K) \leq 2g^*(K).
\]
Now we can run this entire argument for the mirror image $K^!$ giving $s(K^!) \leq 2g^*(K^!)= 2g^*(K)$. Using the properties of $s$ above this implies $-s(K) \leq 2g^*(K)$ so we conclude $|s(K)| \leq 2g^*(K)$ finishing the proof.

\begin{prop}\label{prop:slicetorus}
The slice genus of the $(p,r)$-torus knot is $\frac{(p-1)(r-1)}{2}$.
\end{prop}

Using the $s$-invariant the proof of this is now amazingly simple. It
is clear that the smooth slice genus is less than or equal to the
genus of any Seifert surface. Seifert's algorithm produces a Seifert
surface with Euler characteristic $p-(p-1)r$, that is of genus $(p-1)(r-1)/2$. Thus
\[
|s(T_{p,r})| \leq 2g^*(T_{p,r}) \leq (p-1)(r-1).
\]
But by the exercise above $s(T_{p,r}) = (p-1)(r-1)$ and the result
follows straight away.

The remarkable thing about this proof is that a {\em combinatorially}
defined invariant can tell us something about a result which involves
{\em smoothness}. This is also striking in the following
application on exotic smooth structures.

If you want to prove existence of exotic smooth structure on $\bR^4$
you can do this if you are in possession of a knot which is
topologically slice but not smoothly slice (slice means zero slice
genus). Freedman has a result stating that a knot with Alexander
polynomial 1 is topologically slice. We now have an obstruction ($s$
being non-zero) to being smoothly slice. So armed with these results
all you need to do to calculate $s$ for those knots known to have
Alexander polynomial 1 hoping to reveal one where $s\neq 0$. The knot
in Figure \ref{fig:pretzel}, the $(-3,5,7)$ pretzel knot has $s=-1$.

\begin{figure}[h]
\centerline{
\psfig{figure= 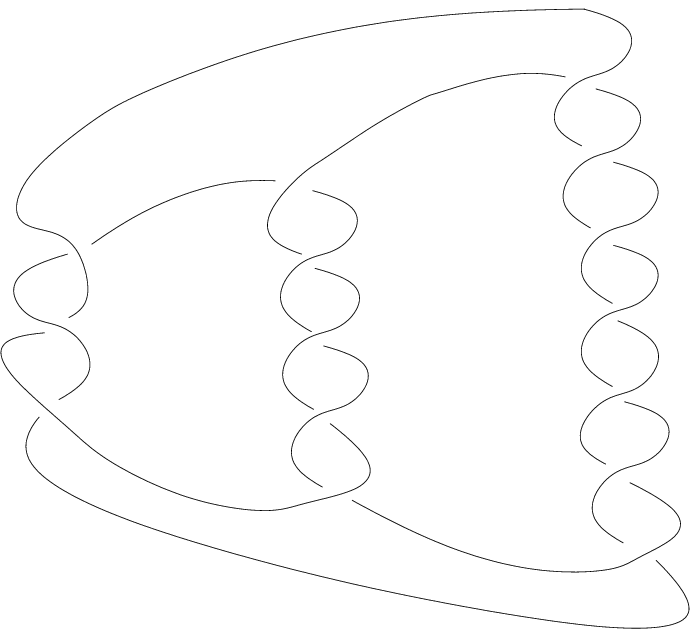}
}
\caption{}
\label{fig:pretzel} 
\end{figure}

\section{Notes and further reading}
\noindent
The family of theories discussed at the beginning of the lecture is
essentially a by-product of Bar-Natan's (geometric) universal theory
\cite{barnatan1}. One applies a TQFT satisfying certain relations to
his theory and thus one only needs to classify the Frobenius algebras
corresponding to these theories. This among other related things can
be found in
\cite{khovanov4}. 
It is possible to work over other rings than $\bQ$. In particular one
can work over the ring $\bZ[h,t]$ to get a universal theory. The
theory defined over $\bZ/2[h]$ (take $t=0$) is known as {\em Bar-Natan
theory}. Proposition \ref{prop:classoverq} is not hard to show - see
\cite{mackaayturnervaz} for details.

\vspace*{2mm}

\noindent
The reference for Lee theory is \cite{lee}

\vspace*{2mm}

\noindent
We have mentioned that Lee theory is filtered rather than bi-graded. An alternative is to work over $\bQ[t]$ where deg$(t)=-4$ and we have $x^2=t1$. This gives a genuine bi-graded theory again. By taking a limit over the ``times $t$'' map one gets back the theory Lee defined. On can compute the bi-graded Lee theory using a spectral sequence (see \cite{turner} for the analogous case of the bi-graded Bar-Natan theory).

\vspace*{2mm}

\noindent
The reference for Rasmussen's invariant is \cite{rasmussen1} where a
proof of Proposition \ref{prop:cobgens} can also be found. Proposition
\ref{prop:slicetorus} was a conjecture (by Milnor) for many years, finally proved in \cite{kronheimermrowka} using gauge theory. Rasmussen's is the first combinatorial proof. In fact there is a more general proposition proved in Rasmussen's paper:  for positive knots the $s$-invariant is
twice the slice genus.

\vspace*{2mm}

\noindent
It was thought for a while (conjectured in \cite{rasmussen1}) that the
$s$-invariant might be equal to twice the $\tau$-invariant in
Heegaard-Floer homology. This is now known to be false and a
counter-example can be found in \cite{heddenording}.


\chapter{Lecture Five}
A natural question is: what else can one categorify? Other knot
polynomials are good candidates. In this lecture we offer a brief
discussion of Khovanov-Rozansky link homology which categorifies a
specialisation of the HOMFLYPT polynomial. We then discuss the topic
of graph homology which has its origins in categorifying graph polynomials. The ``Notes and further reading'' section gives a
cursory look at a number of topics which might be covered in a
hypothetical set of a further five (or more) lectures.

\section{Khovanov-Rozansky homology}
The idea here is to categorify (a specialisation of) the HOMFLYPT
polynomial. The specialisation in question is the one corresponding to
the representation theory of $sl(N)$ and the polynomial
$P_N(D)$ is determined by the skein relation
\[
q^N P_N(\textpicc{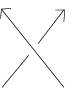})- q^{-N}P_N(\textpicc{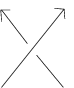}) = (q-q^{-1}) P_N(\textpicc{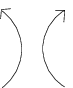}),
\]
and normalised by
\[
 P_N(\text{Unknot}) = \frac{q^N - q^{-N}}{q-q^{-1}}.
\]

To compute the Jones polynomial one can use the Kauffman bracket
which reduces everything to the values (polynomials) assigned to
circles in the plane. For $P_N(D)$ things are not quite as simple,
however one can reduce things to values  assigned to
certain planar graphs. 

Murakami, Ohtsuki and Yamada have defined a polynomial, $P_N(\Gamma)$, for
four-valent planar graphs $\Gamma$ locally modeled on
\textpicc{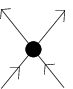}. (In fact they use 3-valent graphs with different
types of edges: elongate the black blob in the picture in the previous
sentence to get a three valent graph looking like
\textpicc{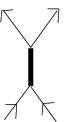}). This polynomial satisfies the following properties:
\begin{enumerate}
\item $P_N(\;\;\textpicc{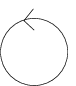}\;\;) = [N ]$ 
\item $P_N(\;\;\textpicc{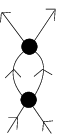}\;\;) = [2] P_N(\;\;\textpicc{blobc.eps}\;\;)$
\item $P_N(\;\;\textpicc{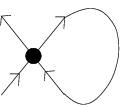}\;\;)= [N-1] P_N(\;\;\textpicc{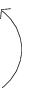}\;\;)$
\item $P_N(\;\;\textpicc{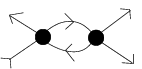}\;\;)= P_N(\;\;\textpicc{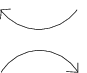}\;\;) + [N-2]P_N(\;\;\textpicc{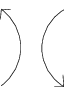}\;\;)$
\item $P_N(\;\;\textpicc{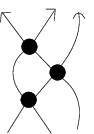}\;\;) + P_N(\;\;\textpicc{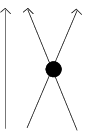}\;\;)=
P_N(\;\;\textpicc{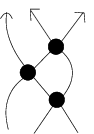}\;\;)+ P_N(\;\;\textpicc{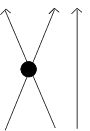}\;\;)$
\end{enumerate}

In the above the square brackets refer to the quantum integer, i.e.  
 \[
 [k] = \frac{q^k - q^{-k}}{q-q^{-1}}.
\]

For an oriented link diagram $D$ resolve each crossing into a 0- or 1-smoothing as indicated in Figure \ref{fig:KRres}.
\begin{figure}[h]
\centerline{
\psfig{figure= 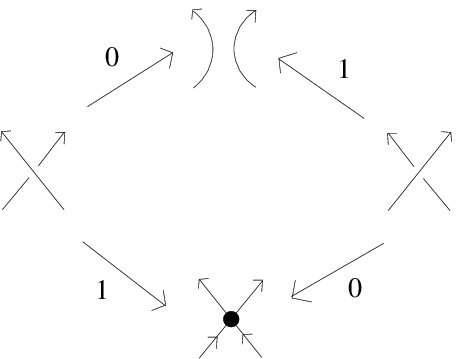}
}
\caption{}
\label{fig:KRres} 
\end{figure}

There are $2^n$ possible smoothings each of which is a planar graph of the sort above. For $\alpha\in\{0,1\}^n$ let $\Gamma_\alpha$ be the associated graph. 
The polynomial $P_N(D)$ can be written as a statesum in terms of the polynomials
$P(\Gamma_\alpha)$ as follows.
\[
P_N(D) = \sum_{\alpha\in\{0,1\}^n}\pm q^{h(\alpha)} P_N(\Gamma_\alpha)
\]
The numbers $h(\alpha)$ and the signs are not hard to determine,
though we will not elaborate on this here.

The categorification of $P_N(D)$ proceeds in two steps: (1) categorify
the polynomial $P(\Gamma)$ i.e. to each graph above assign a vector
space and (2) perform the cube construction to get a polynomial
associated to a link diagram. You need to carry out (1) in such a way
that the (appropriately categorified) properties of $P_N(\Gamma)$ are
satisfied and in such a way that allows you to define maps between two
graphs that differ locally as 0- and 1-smoothings do. You then need
the cube construction to produce a complex whose homology is invariant
under the Reidemeister moves.

It is not surprising that some new ingredients are needed. One such ingredient is the notion of a {\em matrix factorization}. For a commutative ring $R$ and an element $w\in R$, an $(R,w)$-{\em factorization} consists of two free $R$-modules $M^0$ and $M^1$ together with module maps $d^0\colon M^0 \ra M^1$ and $d^1\colon M^1\ra M^0$ such that
\[
d^1\circ d^0 = w Id_{M^0} \spaces \text{ and }\spaces d^0\circ d^1 = w Id_{M^1}.
\]
Put differently, $M=M^0\oplus M^1$ and $d\colon M\ra M$ where
\[
d = \begin{pmatrix} 0 & d^0 \\ d^1 & 0 \end{pmatrix} \spaces \text{ and }\spaces d^2 = wI.
\]
The element $w\in R$ is called the {\em potential}.

\begin{exe}
Take $M=R\oplus R$ and define $d=\begin{pmatrix} 0 & a \\ b & 0 \end{pmatrix}$ for $a,b \in R$. This is an $(R,ab)$-factorization. 
\end{exe}

In fact we want to consider {\em marked} graphs: each arc has one or
more marks on it. Given a marked four-valent planar graph, to each
mark $i$ assign a variable $x_i$ of degree $2$. An example is given in Figure \ref{fig:marks}. 
\begin{figure}[h]
\centerline{
\psfig{figure= 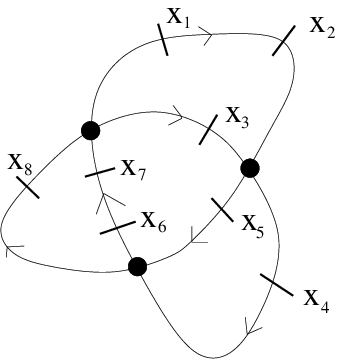}
}
\caption{}
\label{fig:marks} 
\end{figure}

\newpage

We now assign certain matrix factorizations to local pieces of the
graph, which are later tensored together to get something associated
to the graph itself.

We indicate the type of local piece, the ring $R$ and the potential $w$ in the table below.

\vspace{5mm}

\begin{center}
\begin{tabular}{|c|l|l|}
\hline
&&\\
$\;\;\;\;$ Local piece $\;\;\;\;$& $R$ & $w$\\
&&\\
 \hline 
&&\\
\textpicc{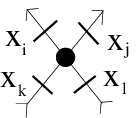} & $\bQ[x_i,x_j,x_k,x_l]$ & $x_i^{N+1} + x_j^{N+1} - x_k^{N+1} - x_l^{N+1}$\\
&&\\
 \hline 
&&\\
\textpicc{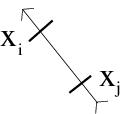} & $\bQ[x_i,x_j]$ &  $x_i^{N+1} - x_j^{N+1}$\\
&&\\
 \hline &&\\
\textpicc{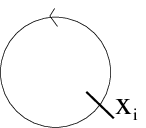} & $\bQ[x_i]$ & $0$\\
&&\\
\hline
\end{tabular}
\end{center}

\vspace{5mm}

It is not particularly enlightening in such a short review to provide
the factorizations explicitly - please refer to the original
source. All these factorizations are then tensored together (over a
variety of intermediate rings) to get a factorization $C(\Gamma)$ which, it turns out,
is a factorization with $R=\bQ[x_i| i \in
\text{set of marks}]$ and $w= 0$. In other words we have 
a length two complex.
The homology of this complex is $\bZ/2 \oplus \bZ$-graded, but is only
non-zero in one of the $\bZ/2$-gradings. We write $H^*(\Gamma)$ for
the homology in this non-zero grading (so $H^*(\Gamma)$ is a 
graded $\bQ$-vector
space). The assignment $\Gamma\mapsto H^*(\Gamma)$ categorifies the
polynomial $P(\Gamma)$.

\begin{prop}
\[
\sum_i q^i \text{dim}(H^i(\Gamma)) = P_N(\Gamma)
\]
\end{prop}

We can now move on to defining a link homology theory. Given an
oriented link diagram $D$ put one mark on each arc of the link and
define 0- and 1-smoothings as in Figure \ref{fig:KRres}. The $2^n$
smoothings $\Gamma_\alpha$ are now marked graphs which as usual we index by
the vertices of the cube $\{0,1\}^n$.

Let $V_\alpha$ be an appropriately shifted version of $H^*(\Gamma_\alpha)$ and set 
\[
C^{i,*} (D) = \bigoplus_{\substack{\alpha\in\{0,1\}^n\\r_\alpha = i+n_+}} V_\alpha .
\]

We need a differential and the key thing is to define the ``partial''
derivatives along edges of the cube. Again we will skip all the details but it is possible to define maps of factorizations as indicated below.
\begin{figure}[h]
\centerline{
\psfig{figure= 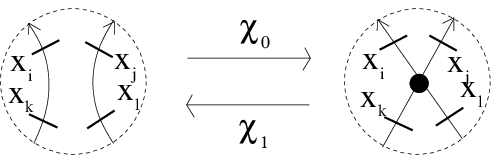}
}
\end{figure}

If $\Gamma_0$ and $\Gamma_1$ are graphs that agree outside a small
region in which they look like the left and right picture above then
$\chi_0$ and $\chi_1$ induce maps $C(\Gamma_0) \ra C(\Gamma_1)$ and
$C(\Gamma_1) \ra C(\Gamma_0)$ and hence maps 
\[
\chi_0\colon H^*(\Gamma_0) \ra
H^*(\Gamma_1) \text{ and } \chi_1\colon H^*(\Gamma_1) \ra H^*(\Gamma_0).
\]
Thus, to each cube edge $\zeta\colon \Gamma \ra \Gamma^\prime$ we can
produce a map $d_\zeta\colon H^*(\Gamma) \ra H^*(\Gamma^\prime)$.
As before set
\[
d  = \sum_{\substack{\zeta \text{ such that}\\ \text{Tail} (\zeta)=\alpha}} \text{sign}(\zeta) d_\zeta.
\]

Miraculously this all works and the following proposition holds.
\begin{prop}
(i) $H(C^{*,*}(D),d)$ is invariant under Reidemeister moves.

(ii) 
\[
\sum_{i,j}(-1)^iq^j \text{dim}(H^{i,j}(D)) = P_N(D).
\]
\end{prop}

Clearly there are many details to check - which is why the paper by Khovanov and Rozansky runs to over one hundred pages!

\section{Graph homology} 
The idea of graph homology is to interpret graph polynomials, like the
chromatic polynomial, Tutte polynomial etc., as the graded Euler
characteristic of a bi-graded vector space. It is much simpler to do
this than it is to work with links, since there is no Reidemeister
invariance to check. None-the-less, graph homology is interesting in
its own right and also serves as a ``toy model'' revealing the same
sort of phenomena that arise in link homology. For example, torsion
also occurs in graph homology and is much more abundant and easier to
get hold of than in link homology.

Let us look at the example of the chromatic polynomial. Let $G$ be a
graph  with vertex set Vert$(G)$ and edge set Edge$(G)$. The {\em
chromatic polynomial}, $P(G)\in \bZ[\lambda]$ is a polynomial which
when evaluated at $\lambda = m \in \bZ$ gives the number of colourings
of the vertices of $G$ by a palette of $m$ colours satisfying the
property that adjacent vertices have different colourings.

There is a procedure to calculate $P(G)$ as follows. Number
 the edges of $G$ by $1,\ldots, n$ and note that there is a one-to-one correspondence between the subsets of edges of $G$ and the set $\{0,1\}^n$. (An edge of $G$ is labelled with 1 if is is present in the subset and 0 otherwise).
For $\alpha\in\{0,1\}^n$ define $G_\alpha$ to be the graph with Vert$(G_\alpha) = \text{Vert}(G)$ and
\[
\text{Edge}(G_\alpha) = \{e_i \in \text{Edge}(G) \mid \text{the $i$'th entry in $\alpha$ is a 1}\}.
\]  

Now define
\[
r_\alpha = \text{ the number of 1's in }\alpha
\]
and 
\[
k_\alpha = \text{ the number of components in }G_\alpha.
\]

A state-sum formula for $P(G)$ is given by
\[
P(G) = \sum_{\alpha\in \{0,1\}^n} (-1)^{r_\alpha} \lambda^{k_\alpha}.
\]

\begin{exercise}
Stop reading here and try to categorify $P(G)$.
\end{exercise}

To categorify $P(G)$ we start with a graded algebra $R$. For
$\alpha\in\{0,1\}^n$ let $R_\alpha = R^{\otimes k_\alpha}$ and as usual
form a cube: associate $R_\alpha$ to the vertex $\alpha$. A simple example is shown in Figure \ref{fig:graphcube}

\begin{figure}[h]
\centerline{
\psfig{figure= 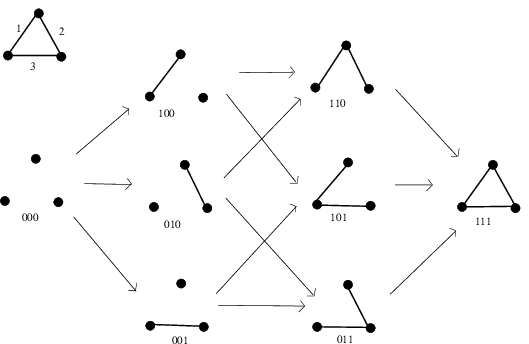}
}
\caption{}
\label{fig:graphcube} 
\end{figure}

Now set
\[
C^{i,*}(G) = \bigoplus_{r_\alpha = i}R_\alpha.
\]
To define a differential we follow the usual procedure. For a cube
edge $\zeta \colon \alpha \ra \alpha^\prime$ note that
$G_{\alpha^\prime}$ either has the same number of components as
$G_\alpha$ or one component less (two components are fused by the
additional edge in $G_{\alpha^\prime}$). Thus we define $d_\zeta\colon
R_\alpha \ra R_{\alpha^\prime}$ to be multiplication in $R$ on  copies of $R$
corresponding to components that fuse (if such exist) and the identity
elsewhere. We thus get a complex whose homology is the graph homology of $G$.

\begin{prop}
\[
\sum_{i,j} (-1)^iq^j\text{dim}(H^{i,j}(G)) = P(G)|_{\lambda \ra \text{qdim}(R)}
\]
\end{prop}

\begin{exercise}
Would taking $d=0$ for the differential work as well?
\end{exercise}

There is a long exact sequence in graph homology which
categorifies the deletion-contraction relation. Given an edge $e$ we
can form two new graphs $G-e$ and $G/e$ where the first has the edge $e$
deleted and the second contracts it. There is a short exact sequence
\[
\xymatrix{
0\ar[r] & C^{i-1,j}(G/e) \ar[r] & C^{i,j}(G) \ar[r] & C^{i,j}(G-e) \ar[r] & 0
}
\]
which gives a long exact sequence

\[
\xymatrix{
\ar[r] & H^{i-1,j}(G/e) \ar[r] & H^{i,j}(G) \ar[r] & H^{i,j}(G-e) \ar[r] & .
}
\]

\section{Notes and further reading} 
\noindent
The first graph polynomial to be categorified was the chromatic
polynomial in \cite{helmeguizonrong}. The dichromatic polynomial was
studied in \cite{stosic}. Torsion in graph homology has been
investigated in \cite{helmeguizonprzytyckirong}.

\vspace*{2mm}

\noindent
The reference for Khovanov-Rozansky theory is
\cite{khovanovrozansky1}. While this theory is considerably harder to
compute than Khovanov's original homology some progress has been
made. In \cite{rasmussen3} Rasmussen describes the Khovanov-Rozansky
polynomial of 2-bridge knots in terms of the HOMFLYPT polynomial and
signature. There is an analogue of Lee's theory investigated by Gornik
in \cite{gornik}. The polynomial of Murakami, Ohtsuki and Yamada is defined and its properties studied in \cite{murakamiohtsukiyamada}.

Khovanov and Rozansky followed up their paper with a sequel
\cite{khovanovrozansky2} in which they consider the two variable HOMFLYPT polynomial. Prior to Khovanov and Rozansky's first paper the case $N=3$ had been
treated in a somewhat different manner by Khovanov in
\cite{khovanov5}. Recently, a link with Hochschild homology has been uncovered \cite{przytycki}.

\vspace*{2mm}

\noindent
Link diagrams can also be drawn on surfaces and the Jones polynomial
can be defined in this context. If the surface $\Sigma$ is part of the
structure then the diagram represents a link in an $I$-bundle over
$\Sigma$. Khovanov homology in this context has been studied in
\cite{asaedaprzytyckisikora}. If the surface is not really part of the structure, but rather just a carrier for the diagram (so you can add/subtract handles away from the diagram) then equivalence classes of diagrams are known as {\em virtual links}. Khovanov homology of these has been studied in \cite{manturov} and \cite{turaevturner}.

\vspace*{2mm}

\noindent
The Jones polynomial corresponds to the 2-dimensional representation
of $U_q(sl_2)$ and allowing other representations leads to the {\em
coloured Jones polynomial}. A link homology theory categorifying this
was defined in \cite{khovanov3}.

\vspace*{2mm}

\noindent
One of the most interesting questions surrounding the subject is to
uncover the geometry that lies behind Khovanov homology. A proposal for a framework unifying Khovanov-Rozansky homology and knot Floer homology has be put forward in \cite{dunfieldgukovrasmussen}.

In a different direction P. Seidel and I. Smith 
have constructed a homology theory for links using symplectic geometry
\cite{seidelsmith}. This theory is conjectured to be isomorphic to
Khovanov homology (after suitably collapsing the bi-grading into a
single grading). Building on this Manolescu has constructed a similar theory
for each $N$ and conjectured this to be isomophic to Khovanov-Rozansky
homology \cite{manolescu}.

\vspace*{2mm}

\noindent
Another exciting direction is to try to give some ``physical''
interpretation for Khovanov homology (in the sense that Witten gave a
physical interpretation of the Jones polynomial as the partition
function of a quantum field theory). S. Gukov, A. Schwartz and C. Vafa
have made an attempt in this direction \cite{gukovschwartzvafa}
conjecturing a connection to string theory.

\vspace*{2mm}

\noindent
Last, but certainly not least, there has been a huge effort to write
computer programs to calculate Khovanov homology groups. The first of
these by Bar-Natan (using Mathematica) coped with links up to 11 or 12
crossings. This was improved on by Shumakovitch \cite{shumakovitchKhoHo} 
with a program using
Pari. Bar-Natan now has a nice theoretical trick which
speeds things up considerably. This has been implemented by Jeremy Green 
and you can download the package at from
the (wonderful) knot atlas (set up by Dror Bar-Natan and Scott
Morrison).
\begin{center}
http://katlas.math.toronto.edu/wiki/
\end{center}


\end{document}